\documentclass{article}
\usepackage{amssymb,amsmath}
\usepackage{graphics}
\input{epsf}
\input xy

\xyoption{all}
\usepackage[all,cmtip]{xy}
\usepackage{graphicx}
\usepackage{amsmath,amsthm,amsfonts,amssymb}
\usepackage{epstopdf}  

\usepackage{graphics}
 \input xypic
 \input epsf
 \input xy
 \xyoption{all}

\begin{document}

\def\adots{\mathinner{\mkern2mu\raise1pt\hbox{.}\mkern2mu\raise4pt\hbox{.}
\mkern2mu\raise7pt\hbox{.}\mkern1mu}}

\centerline{\huge Charles University in Prague}

\bigskip

\centerline{\huge Faculty of Mathematics and Physics}

\bigskip

\bigskip

\bigskip

\bigskip

\bigskip

\bigskip

\bigskip

\bigskip

\bigskip

\centerline{\mbox{\includegraphics[width=60mm]{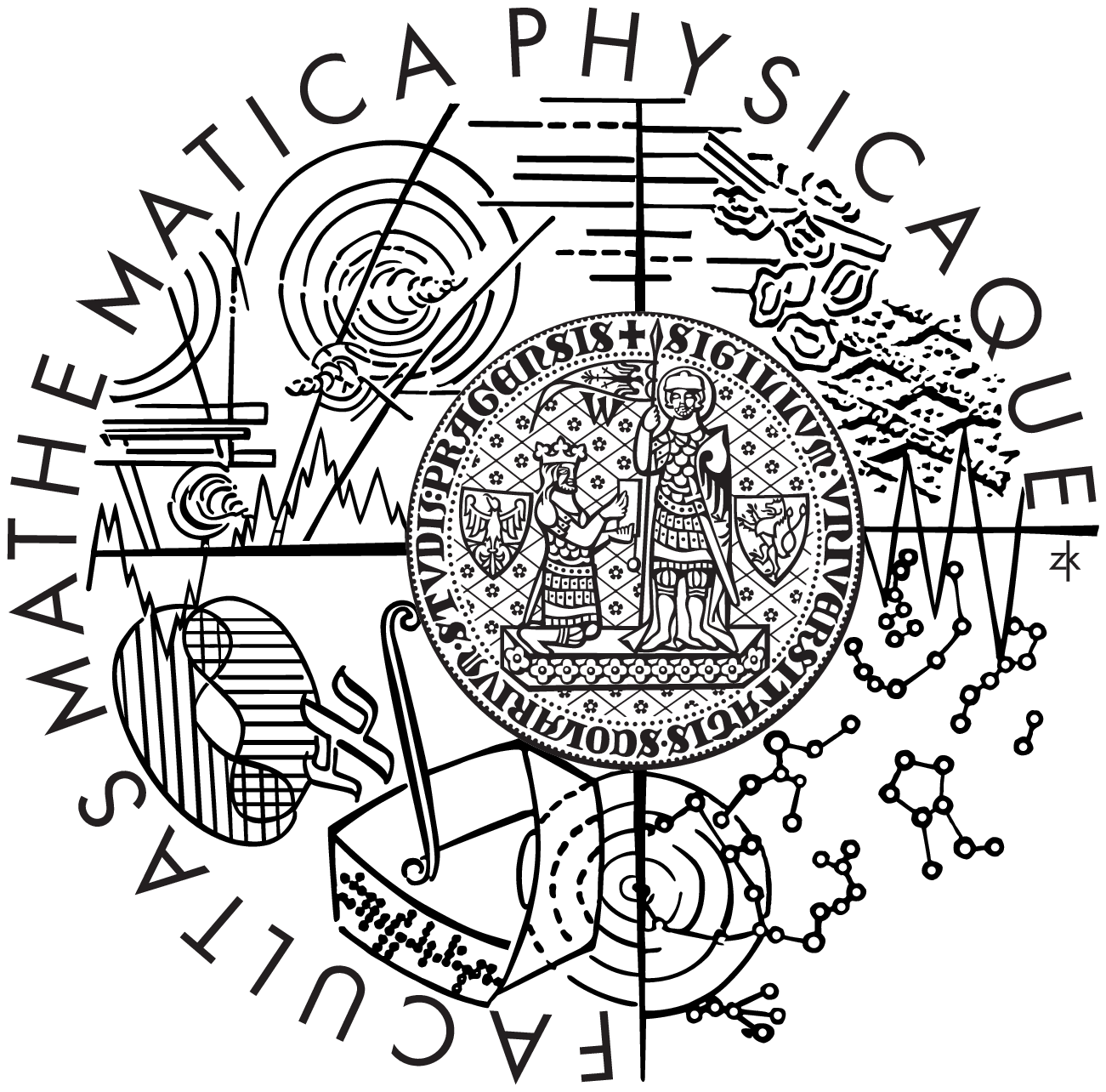}}}

\bigskip

\bigskip

\bigskip

\bigskip

\bigskip

\bigskip

\bigskip

\bigskip

\bigskip

\bigskip

\centerline{\huge  Habilitation Thesis}

\bigskip

\centerline{\huge \bf Symplectic spinors and  Hodge theory}

\bigskip

\bigskip

\bigskip

\bigskip

\centerline{\huge Svatopluk Kr\'ysl}


 \bigskip
 
 \bigskip
 
 \bigskip

\bigskip

\bigskip

\centerline{\huge Discipline: Mathematics  - Geometry and Topology} 

\bigskip

\centerline{\huge June, 2016}
 
\thispagestyle{empty}

\newpage

\thispagestyle{empty}

\phantom{S} 

\newpage

\thispagestyle{empty}


\newpage

\tableofcontents

\newpage

\setcounter{page}{1}

\section{Introduction}

 \bigskip 
 
 \bigskip

 {\it I have said the 21st century might be the era of
quantum mathematics or, if you like, of infinite-dimensional mathematics. What
could this mean? Quantum mathematics could mean, if we get that far, `understanding
properly the analysis, geometry, topology, algebra of various non-linear
function spaces', and by `understanding properly' I mean understanding it in such
a way as to get quite rigorous proofs of all the beautiful things the physicists have
been speculating about.}

\vspace{0.6cm}

{\hspace{11.2cm} \it  Sir Michael Atiyah} 
  
\bigskip

  In the literature, topics contained in this thesis are treated rather separately. 
From a philosophical point of view, a common thread of themes that we consider  is represented by the above quotation of M. Atiyah.
We are inspired by mathematical and theoretical physics.

\bigskip

  The content of the thesis concerns  analysis, differential geometry  and representation theory on infinite dimensional objects.
A specific infinite dimensional object  which we consider is the Segal--Shale--Weil  representation of the metaplectic group. This representation  
originated in number theory and theoretical physics.  
We analyze its tensor products with finite dimensional representations, 
induce it to metaplectic structures defined over
symplectic and  contact projective manifolds  obtaining differential operators whose properties we study.
From other point of view,  differential geometry uses the infinite dimensional algebraic objects to obtain vector bundles and differential operators, 
that we investigate by generalizing analytic methods known currently -- 
namely by a Hodge theory for complexes in categories of specific modules over $C^*$-algebras. 
\bigskip

Results described in the thesis arose from the year 2003 to the year 2016. 
At the beginning, we aimed to classify all {\bf first order invariant differential operators} acting between bundles over contact projective manifolds that are induced by   
those irreducible representations  of the metaplectic group which have bounded 
 multiplicities. See Kr\'ysl \cite{KryslDGA1} for a result. 
Similar results were achieved by Fegan \cite{Fegan} in the case of irreducible finite rank bundles over 
Riemannian manifolds equipped with a conformal structure. In both cases, for any such two  bundles, there is at most one first order invariant 
differential operator up to a scalar multiple.\footnote{and up to operators of the zeroth order. See section 4.3.} The condition for the existence in the case of contact projective manifolds is expressed by the highest weight of the induced representation
considered as a module over a suitable simple group, by a conformal weight, 
and by the action of $-1 \in \mathbb{R}^{\times}.$ The result  is based on a 
decomposition of the tensor product of an irreducible highest weight 
$\mathfrak{sp}(2n,\mathbb{C})$-module that has bounded multiplicities with the 
defining representation of $\mathfrak{sp}(2n,\mathbb{C})$ into 
irreducible submodules. See Kr\'ysl \cite{KryslJOLT1}. 
For similar classification results in the case of more general  parabolic geometries and bundles induced by finite dimensional modules, see Slov\'ak, Sou\v{c}ek \cite{SlovakSoucek}.

\bigskip
  
Our next research interest, described in the thesis, was the {\bf de Rham sequence tensored by the Segal--Shale--Weil representation}.  
From the algebraic point of view, the Segal--Shale--Weil representation (SSW representation)  is  an 
$L^2$-globalization of the direct sum of two specific infinite dimensional 
Harish-Chandra $(\mathfrak{g}, K)$-modules with bounded multiplicities over the metaplectic Lie algebra, which are called completely pointed.
 We decompose the  de Rham 
sequence with values in the mentioned direct sum into sections of irreducible bundles, i.e., bundles induced by 
irreducible representations. See Kr\'ysl \cite{KryslArchivumSSVF}. For a $2n$ 
dimensional symplectic manifold with a metaplectic structure, 
the bundle of exterior forms of degree $i$, $i \leq n,$ with values in the 
Segal--Shale--Weil representation decomposes into $2(i+1)$ irreducible bundles. 
For $i \geq n,$ the number of such bundles is $2(2n-i+1).$
It is known that the decomposition structure of completely reducible 
representations is  connected to the so-called Schur--Weyl--Howe dualities. See 
Howe \cite{HoweInvariant} and Goodman, Wallach \cite{GoodmanWallach}.
We use the decomposition of the twisted de Rham sequence to
obtain a duality for the metaplectic group which acts in this case, on the 
exterior forms with values in the  SSW representation.
See Kr\'ysl \cite{KryslJOLT2}.  The dual partner to the metaplectic group appears to be the orthosymplectic Lie superalgebra $\mathfrak{osp}(1|2).$ 

\bigskip

Any Fedosov connection (i.e., a symplectic and torsion-free connection) on a symplectic mani\--fold with a metaplectic structure 
defines a covariant derivative on the symplectic spinor  bundle which is the bundle induced by the Segal--Shale--Weil representation. We prove that twisted de Rham derivatives  
 map sections of 
an irreducible subbundle into sections of at most three irreducible subbundles. 
Next, we are interested in the quite fundamental question on the 
structure of the curvature tensor of the symplectic spinor covariant derivative 
similarly  as one does in the study of the curvature of a Levi-Civita or a Riemannian connection. 
Inspired by  results of Vaisman in \cite{Vaisman} on the curvature tensors of Fedosov connections, we derive a decomposition of the
curvature tensor on  symplectic spinors. See Kr\'ysl \cite{KryslJGP1}.
Generalizing this decomposition, we are able to find certain subcomplexes  of the twisted de Rham sequence, that 
are called symplectic twistor complexes in a parallel to the spin geometry.
These complexes exist under specific restrictions on the curvature of the Fedosov connection. Namely, the connection is demanded to be 
of  Ricci-type.
See Kr\'ysl \cite{KryslMonat}.
Further results based on the decomposition of the curvature concern  a
relation of the spectrum of the symplectic spinor Dirac operator to the spectrum of the symplectic spinor Rarita--Schwinger operator.
See Kr\'ysl \cite{KryslArchivumRSDR}.   
The symplectic Dirac operator was introduced by K. Habermann. See 
\cite{HabermannAGAG}. 
  The next result is on symplectic  Killing spinors. We 
prove that if there exists a non-trivial  (i.e., not  covariantly constant) symplectic 
Killing spinor, the connection is not Ricci-flat.
See \cite{KryslCMUC}.

\bigskip

Since the 
classical theories on analysis of elliptic operators  on compact manifolds 
are not applicable in the case of the  de Rham complex twisted by the Segal--Shale--Weil representation,
we tried  to develop a {\bf Hodge theory for infinite 
rank bundles}. 
We use and elaborate ideas of Mishchenko and  Fomenko (\cite{Mishchenko} and \cite{MishchenkoFomenko}) on generalizations 
of the Atiyah--Singer index theorem to investigate cohomology groups
of infinite rank elliptic complexes  concerning their topological and algebraic properties.
We work in the categories $PH_A^*$ and $H_A^*$ whose  objects  are pre-Hilbert 
$C^*$-modules and Hilbert $C^*$-modules, respectively,  and whose morphisms are adjointable maps 
between the objects.
These notions go back to the works of Kaplansky \cite{Kaplansky}, Paschke \cite{Paschke} and Rieffel \cite{Rieffel}. 
Analyzing proofs of the classical Hodge theory, we are lead to the notion of 
a Hodge-type complex in an  additive and dagger category.
We introduce a class of  self-adjoint parametrix possessing 
complexes, and prove that any self-adjoint parametrix possessing complex in $PH_A^*$ is of 
Hodge-type.
Further, we prove that in $H_A^*$  the category of self-adjoint 
parametrix possessing complexes coincides with
the category of the Hodge-type ones. 
Using these results, we show that
an elliptic complex on sections of finitely generated projective Hilbert $C^*$-bundles over compact manifolds are of Hodge-type if the images of 
the Laplace operators of the complex are closed. 
The cohomology groups of such complexes are isomorphic to the kernels of  
the Laplacians and they 
are Banach spaces with respect to the quotient topology.  
Further, we prove 
that the cohomology groups are finitely generated projective Hilbert $C^*$-modules. See Kr\'ysl \cite{KryslAGAG1}.
Using the result of Baki\'c and Gulja\v{s} \cite{BakicGuljas} for modules over a 
$C^*$-algebra of compact operators $K$, we are able to remove the condition on 
the closed image. We prove in \cite{KryslJGP2} that elliptic 
complexes of differential operators on finitely generated projective 
$K$-Hilbert bundles are of 
Hodge-type and that their cohomology groups are finitely generated projective 
Hilbert $K$-modules.   See Kr\'ysl  \cite{KryslDGA2}, \cite{KryslJGSP} and \cite{KryslAGAG2} for a possible application 
connected to the quotation of Atiyah. 

\bigskip

We find it more appropriate to mention author's results and their context
in Introduction, and treat motivations,  development of
important notions, and  most of the references in  Chapters 2 and 3.  
In the 2nd  Chapter, we recall a definition,   realization and characterization of the Segal--Shale--Weil representation.
In Chapter 3, symplectic manifolds,  Fedosov connections,  metaplectic structures, symplectic Dirac and certain related operators are introduced.
Results of K. and L. Habermann on global analysis related to these operators are mentioned in this part as well.
Chapter 4 of the thesis contains  own results  of the 
applicant. We start with the appropriate representation theory and Howe-type duality.
Then we treat results on the twisted de Rham derivatives, curvature of the symplectic spinor derivative and twistor complexes. Symplectic Killing spinors are defined in this part as well. 
We give a classification of 
the invariant operators for contact projective geometries  together with   results 
on the decomposition of the appropriate tensor products in the third subsection.
The fourth subsection  is devoted to the Hodge theory.   
The last part of the thesis consists of selected articles published in 
the period 2003--2016.

\newpage

\section{Symplectic spinors}

  The discovery of   symplectic spinors as a rigorous mathematical object is attributed to
I. E. Segal, D. Shale and A. Weil. See Shale \cite{Shale} and  Weil \cite{Weil}. 
Segal and Shale considered the real symplectic group as a group of 
canonical transformations and were interested in a certain quantization of Klein--Gordon fields. 
Weil was interested in number theory connected to
 theta functions, so that he considered more general symplectic groups than the ones over the real numbers. Let us notice, that this representation appeared in 
works of van Hove (see Folland \cite{Folland}, p. 170) at the Lie algebra 
level and can be found in certain formulas of Fresnel in wave optics  already 
(see Guillemin, Sternberg \cite{GuilleminSternberg}, p. 71).

\bigskip

  When doing quantization, one has to
 assign to ``any'' function defined
 on the phase space of a considered classical system
an operator acting on a certain function space -- a Hilbert space by a rule. Usually, 
smooth functions are considered to represent the right  class for the set of  quantized functions.
The prescription shall assign to the Poisson bracket
of two smooth functions a multiple of the commutator of the
operators assigned to the individual functions. The 
multiple is determined by ``laws of nature''. 
It equals to $(\imath\hslash)^{-1},$ where $\imath$ is a fixed root of $-1$ and 
$\hslash$ is the Planck constant over
$2\pi$.  Thus, in the first steps, the quantization map is demanded to 
be a Lie algebra homomorphism up to a multiple. In further considerations, there 
is a freedom allowed in the sense that the image of a Poisson bracket  need not be 
the $(\imath\hslash)^{-1}$ multiple of the commutator precisely, but the so-called 
deformations are allowed. (See Waldmann 
\cite{Waldmann} and also Markl, Stasheff \cite{MarklStasheff} for a framework of quite general deformations.)
This tolerance is mainly because of the 
Groenewold--van Hove no go theorem (see Waldmann \cite{Waldmann}).    Analytically, 
deformations
are convergent series in the small variable $\hslash.$ Their connection to 
the formal deformations is given by the so called Borel lemma \cite{Waldmann}.

\bigskip

 The state space of a classical system with finite degrees of freedom is modeled 
by a symplectic manifold $(M,\omega).$
The state space of a point particle moving in an $n$-dimensional vector space or $n$
 point particles on a line is the space
$\mathbb{R}^{2n}$ or the intersection of open half-spaces 
in it, respectively. Considering the coordinates
$q^1,\ldots, q^n,$ and $p_1, \ldots, p_n$ on $\mathbb{R}^{2n}$ where $q^i$ projects onto the $i$-th 
coordinate and $p_i$ onto the $(n+i)$-th one,  $\omega$ 
equals to $\sum_{i=1}^n dq^i \wedge dp_i,$ or to its restriction to the 
intersection, respectively.

\bigskip

The  group of all linear maps $\Phi$ on $\mathbb{R}^{2n}$ which preserve the 
symplectic form, i.e., $\Phi^*\omega  = \omega,$ is called the {\bf symplectic group}.
Elements of this group do not change the form of dynamic equations 
governing  non-quantum systems -- the Hamilton's equations. In  this 
way, they coincide with linear canonical transformations used in Physics.\footnote{The system is supposed to be non-dissipative, i.e., its 
Hamiltonian function does not depend on time.} See, e.g., Arnold \cite{Arnold} or Marsden, Ratiu 
\cite{MarsdenRatiu}.

\bigskip

The symplectic group $G=Sp(2n,\mathbb{R})$ is an $n(2n+1)$ dimensional Lie 
group.  Its maximal compact subgroup is isomorphic to the unitary 
group $U(n)$ determined by choosing a 
compatible positive linear complex structure, i.e., an $\mathbb{R}$-linear map $J: \mathbb{R}^{2n} \to \mathbb{R}^{2n}$ such that 1)
$J^2 = -1_{\mathbb{R}^{2n}}$ and 2) the bilinear map $g: \mathbb{R}^{2n} \times \mathbb{R}^{2n} \to \mathbb{R}$ given by
$g(v,w)=\omega(v,Jw)$ is symmetric and positive definite, i.e., a scalar product.  
The unitary group  can be proved  diffeotopic to the circle $S^1,$ and consequently, its first fundamental group is isomorphic to $\mathbb{Z}.$ Thus, 
for each $m\in \mathbb{N},$ $Sp(2n,\mathbb{R})$ has a 
unique non-branched $m$-folded covering  
$\lambda(m): {\widetilde{\mbox{}^m Sp(2n,\mathbb{R})}} \to Sp(2n,\mathbb{R})$ up to a covering isomorphism.
We fix an element $e$ in the preimage of the neutral element in $Sp(2n, \mathbb{R})$ on the two fold covering.  The unique Lie group  
such that its neutral element is $e$ and such that the 
covering map is a Lie group homomorphism is  called the {\bf metaplectic 
group}, or if we wish, the real metaplectic group.   We set $\lambda = \lambda(2)$ and 
$\widetilde{G}=Mp(2n,\mathbb{R}) = \widetilde{\mbox{}^2\,Sp(2n,\mathbb{R})}.$  
We denote the $\lambda$-preimage of $U(n)$ by $\widetilde{K}.$

\bigskip

So far, the construction of the metaplectic group was rather abstract. One of 
the basic results of the structure theory of Lie groups  is
that this is unavoidable indeed. By this we mean that there is no 
faithful representation of $Mp(2n,\mathbb{R})$ by matrices on 
a {\bf finite dimensional} vector space. Otherwise said, the metaplectic group cannot be realized as 
a Lie subgroup of any finite dimensional general linear group. 
Following Knapp \cite{KnappBeyond}, we prove this statement.  

\bigskip

\noindent{\bf Theorem 1:} The connected double cover $Mp(2n,\mathbb{R})$ does not have a 
realization as a Lie subgroup of $GL(W)$ for a finite dimensional vector space $W.$

\noindent{\it Proof.} Let us suppose that there exists a faithful representation $\eta': \widetilde{G} \to \mbox{Aut}(W)$
of the metaplectic group on a finite dimensional space $W.$ This 
representation gives rise to a faithful representation
$\eta: \widetilde{G} \to \mbox{Aut}(W^\mathbb{C})$ on the complexified vector 
space $W^{\mathbb{C}}.$ This map is  injective by its construction.
The corresponding Lie algebra representation, i.e., the map 
$\eta_*:\widetilde{\mathfrak{g}} \to \mbox{End}(W^{\mathbb{C}})$ is well defined because of the 
finite dimension of $W^{\mathbb{C}}.$
Consequently, we have the commutative diagram
$$\begin{xy}\xymatrix{
\widetilde{\mathfrak{g}} \ar[d]^-{\lambda_*} \ar[r]^-{\eta_*}& \mbox{End}(W^\mathbb{C})\\
\mathfrak{g} \ar[r]^-{j'}&  \mathfrak{sp}(2n,\mathbb{C})\ar[u]^{\phi'}
}\end{xy}$$
where $\mathfrak{g}$ is the Lie algebra of the appropriate symplectic group, 
$j'$ is the natural inclusion and $\phi'$ is defined by $\phi'(A + \imath B)= \eta_*\lambda_*^{-1}(A) + \imath \eta_*\lambda^{-1}_*(B),$ $A,B\in \mathfrak{g}.$ 
Taking the exponential of the Lie algebra $\mathfrak{sp}(2n,\mathbb{C}) 
\subseteq \mbox{End}(\mathbb{C}^{2n}),$
we get the  group $Sp(2n,\mathbb{C}).$ Because this group is simply connected, we get a 
representation
$\phi: Sp(2n,\mathbb{C}) \to \mbox{Aut}(W^{\mathbb{C}})$ which integrates $\phi'$ in the sense that
$\phi_*=\phi'.$ 
Because $\lambda_*,$ $\eta_*$ and also $\phi_*$ 
are derivatives at $1$ of the corresponding Lie groups representations, and $j'$ is the derivative at $1\in G$
of the canonical inclusion $j: G \to Sp(2n,\mathbb{C}),$ 
we obtain a corresponding diagram at  Lie groups  level which is commutative as well.
Especially, we have $\eta = \phi \circ j \circ \lambda.$ However,  the right 
hand side of this expression is not injective whereas the left hand side is.    
This is a contradiction. \hfill $\Box$

\bigskip
  
\noindent{\bf Remark:} The complex orthogonal group $SO(n,\mathbb{C})$ is not simply connected,
so that the above proof does not apply for $G=SO(n,\mathbb{R})$ and its connected double cover
$\widetilde{G} = Spin(n,\mathbb{R}).$  If it applied,  $Spin(n,\mathbb{R})$ would 
not have any faithful finite dimensional representation.

\subsection{The Segal--Shale--Weil representation}

For the canonical symplectic vector space $(\mathbb{R}^{2n},\omega),$ a group 
$H(n)=\mathbb{R}^{2n} \times \mathbb{R}$ is assigned in which the
group law is given by
$$(v,t) \cdot (w,s) = (v + w, s + t + \frac{1}{2}\omega(v,w))$$ 
where $(v,t), (w,s) \in H(n).$ This is the so called 
{\bf Heisenberg group} $H(n)$ of dimension $2n+1.$
Let us set $L=\mathbb{R}^n \times \{0\} \times \{0\} \subseteq H(n)$  and $L'=\{0\}\times \mathbb{R}^n \times \{0\}$ for the vector space of 
initial space  and  for the vector space of 
initial impulse conditions, respectively. In particular,
the symplectic vector space $\mathbb{R}^{2n}$ is isomorphic to the direct sum 
$L\oplus L'.$

\bigskip
 
For any $u'\in H(n),$ we have a unique $t\in \mathbb{R}$ and $q \in L,$  $p 
\in L'$ such that $u' = (q, p, t).$
The {\bf Schr\"odinger representation} $Sch$ of the Heisenberg group
$Sch: H(n) \to U(L^2(L))$ is given by $$(Sch((q,p, t))f)(x) = e^{2\pi \imath t + \pi \imath \omega(q,p) + 2 \pi \imath \omega(x,p)}f(x+q)$$ 
where $q, x \in L,$ $p\in L',$ $t\in \mathbb{R},$ and $f\in L^2(L).$ It is an irreducible representation. See Folland \cite{Folland}.
(By $\omega(x,p),$ we mean $\omega((x,0),(0,p))$ and similarly for $\omega(q,p).$)

We may ``twist'' this representation in the following way. For any $g\in G,$ we set
$l_g: H(n) \to H(n),$ $l_g(u,t)=(gu,t),$ where $u\in \mathbb{R}^{2n}$ and $t\in \mathbb{R}.$
Consequently,  we obtain a family of representations $Sch \circ l_g$ of the Heisenberg group $H(n)$  
parametrized by elements  $g$ of the symplectic group $G.$ 
The action of the center of $H(n)$ is the same for each member of the  family $(Sch \circ l_g)_{g\in G}.$ 
Namely, $(Sch \circ l_g)(0,0,t)= e^{2 \pi \imath t},$ $t \in \mathbb{R}.$  
Let us recall the Stone--von Neumann theorem. For a proof, we refer to Folland \cite{Folland} or Wallach \cite{Wallach}.

\bigskip

\noindent{\bf Theorem 2} (Stone--von Neumann): Let $T$ be an irreducible 
unitary representation of the Heisenberg group on a separable infinite 
dimensional Hilbert space $H.$   Then $T$ is unitarily 
equivalent to the Schr\"odinger representation.

\bigskip

By Stone--von Neumann theorem, we find a unitary 
operator $C_g: L^2(L) \to L^2(L)$  that  intertwines  $Sch\circ l_g$  and $Sch$ for each $g\in G.$\footnote{We say that
$C: W \to W$ intertwines a representation $T: H \to \mbox{Aut}(W)$ of the group $H$ if $C\circ T(h)= T(h) \circ C$ for each  $h \in H.$}
 By Schur lemma for irreducible 
unitary representations (see Knapp \cite{KnappRed}), we see that there is a function 
$m: G \times G \to U(1)$ such that  $C_g C_{g'} = m(g,g')C_{gg'},$ $g,g'\in G.$
In particular, $g \mapsto C_g$ is a projective representation of
$Sp(2n,\mathbb{R})$ on the Hilbert space $L^2(L).$
It was proved by Shale in \cite{Shale} and Weil in \cite{Weil} that it is 
possible to lift the cocycle $m$  and the projective representation $g \mapsto C_g$ 
of $G$ to the metaplectic group to obtain a true representation of the $2$-fold cover.  
We denote this representation by $\sigma$ and call it the {\bf Segal--Shale--Weil representation}.
Note that some authors call it the symplectic spinor,   
metaplectic or  oscillator representation.  The representation is unitary and
faithful.  See, e.g., Weil \cite{Weil}, Borel, Wallach \cite{BorelWallach}, Folland \cite{Folland}, Moeglin et 
al.  \cite{MoeglinVignerasWaldspurger},
Habermann, Habermann \cite{HabermannHabermann} or Howe \cite{Howetheta}.
 
\bigskip

The ``essential'' uniqueness of the Segal--Shale--Weil representation with respect to the choice of a representation of the Heisenberg group
is expressed in the next theorem.

\bigskip

\noindent{\bf Theorem 3:} Let $T: H(n) \to U(W)$ be an irreducible  unitary representation
of the Heisenberg group on a Hilbert space $W$ and $\sigma': \widetilde{G} \to U(W)$ be a non-trivial unitary representation
of the metaplectic group such that for all $(v,t)\in H(n)$ and $g\in \widetilde{G}$ 
$$\sigma'(g) T(v,t)\sigma'(g)^{-1} = T(\lambda(g)v,t).$$ Then there exists a deck transformation $\gamma$
of $\lambda,$  such that $\sigma'$ is equivalent either to
$\sigma \circ \gamma$ or to $\sigma^* \circ \gamma,$ where $\sigma^*(g)=\tau \sigma(g) \tau$ and
$(\tau(f))(x) = \overline{f(x)},$ $x \in \mathbb{R}^n,$ $g\in \widetilde{G}$ and $f\in L^2(\mathbb{R}^n).$ 

\noindent{\it Proof.}  See Wallach \cite{Wallach}, p. 224. \hfill $\Box$

\bigskip

\noindent{\bf Remark:} Let us recall that a deck transformation $\gamma$ is any continuous map 
which satisfies  $\lambda \circ \gamma = \lambda.$ Note that in the case of the symplectic group covered by the 
metaplectic group, a deck transformation is either the identity map or 
the map ``interchanging'' the leaves of the metaplectic group.

\subsection{Realization of symplectic spinors}

There are several different objects that one could call a symplectic basis.
We choose the one which is convenient for considerations in projective contact geometry. 
(See Yamaguchi \cite{Yamaguchi} for a similar choice.) 
If $(V,\omega)$ is a symplectic vector space of dimension $2n$ over a field 
$\Bbbk$ 
of characteristic zero, we call a basis $(e_i)_{i=1}^{2n}$ of $V$ a {\bf 
symplectic 
basis} if
$\omega(e_i,e_j)=\delta_{i, 2n+1-j}$ for $1\leq i \leq n$ and $1 \leq j \leq 
2n,$ and $\omega(e_i,e_j)=-\delta_{i,2n+1-j}$ for $n+1 \leq i \leq 2n$ and $1\leq j \leq 2n.$
Thus, with respect to a symplectic basis, the matrix of the symplectic form is
{\small{$$
(\omega_{ij})=\left(\begin{array}{c|c}
  0 & K \\ 
  \hline
  -K & 0
  \end{array}\right)
$$} }

where $K$ is the following $n\times n$ matrix 

$$K= \left({\begin{array}{cccc}
      0 & \ldots & 0      & 1\\
      0 & \ldots & 1     & 0\\
      \vdots &  \adots & \vdots&\vdots\\
      1 &   \ldots &0 & 0
             \end{array}} \right).$$
             
\bigskip

Further, we denote by  $\omega^{ij},$ $i,j=1,\ldots, 2n,$ the coordinates which 
satisfy
$\sum_{k=1}^{2n}\omega_{ik}\omega^{jk}= \delta_i^j$ for each $i,j=1,\ldots,2n.$ 
They define a bilinear form 
$\omega^*: V^* \times V^* \to \Bbbk,$ e.g., by setting $\omega^*= 
\sum_{i,j=1}^{2n}\omega^{ij} e_i \wedge e_j.$
We use $\omega_{ij}$ and $\omega^{ij}$ to {\bf rise} and {\bf lower} indices of 
tensors over $V.$
 For coordinates ${K_{ab\ldots c\ldots d}}^{rs \ldots t \ldots u}$ of a tensor 
$K$ on $V,$  we denote
the expression $\sum_{c=1}^{2n}\omega^{ic}{K_{ab\ldots c \ldots d}}^{rs \ldots t}$ by 
${{{K_{ab \ldots}}^{i}}_{\ldots d}}^{rs \ldots t}$ and 
$\sum_{t=1}^{2n}{K_{ab\ldots c}}^{rs \ldots t \ldots u}\omega_{ti}$ by ${{{K_{ab \ldots 
c}}^{rs\ldots}}_{i}}^{\ldots u}$ and similarly for other types of
tensors and  in the geometric setting when we consider tensor fields on 
symplectic manifolds.

\bigskip
             
\noindent{\bf Remark:} Let $(\mathbb{R}^{2n},\omega)$ be the canonical 
symplectic vector space introduced at the beginning of this Chapter.
Then the canonical arithmetic basis of $\mathbb{R}^{2n}$ is not a symplectic 
basis according to our definition unless $n=1.$

\bigskip
 
 Let us denote  the $\lambda$-preimage of $g\in 
Sp(2n,\mathbb{R})$ by $\widetilde{g}.$
Suppose $A, B \in M_n(\mathbb{R}),$  $A$ is invertible and $B^t=B.$ 
 We define the  following 
representation
of $\widetilde{G}$ on $L^2(\mathbb{R}^n)$ 
\begin{eqnarray*}
(\sigma(h_1)f)(x) &=& \pm e^{- \pi \imath g_0(Bx,x)/2}f(x) 
\mbox{ for any } h_1 \in \widetilde{g_1}, \,   g_1=\left( 
\begin{array}{c|c}
1 & 0 \\
\hline
B & 1
\end{array} \right)\\
(\sigma(h_2)f)(x) &=&  \sqrt{\mbox{det} \, A^{-1}}f(A^{-1}x)
\mbox{ for any } h_2 \in \widetilde{g_2}, \, 
g_2= \left(\begin{array}{c|c}
A  &  0 \\

\hline

0 & {A}^{-1t}
\end{array} \right) \\
  (\sigma(h_3)f)(x)&=& \pm \imath^n e^{\pi \imath n/4} (\mathcal{F}f)(x) \mbox{ for any } 
h_3 \in \widetilde{g_3}, \,  g_3 = J_0 =\left(
\begin{array}{c|c}
  0 & -1 \\ 
  \hline
  1 & 0
  \end{array}\right)
\end{eqnarray*}	 
where  $f \in L^2(\mathbb{R}^n)$  and $x\in \mathbb{R}^n.$ The $\pm$ signs and the square roots in the 
definition
of $\sigma(h_i)$  depend on the specific element in the 
preimage of $g_i.$ The coordinates   of $g_i,$ 
$i=1,2,3,$  are considered
with respect to the canonical basis of   $\mathbb{R}^{2n}.$  See 
Folland \cite{Folland}.
Notice that we use the Fourier transform  defined by 
$(\mathcal{F}f)(y) =\int_{x\in \mathbb{R}^n} e^{-2\pi\imath g_0(x, y)} f(x) \mathrm{d} x,$  $y\in \mathbb{R}^n,$ 
with respect to the
Lebesgue measure $\mathrm{d} x$ on $\mathbb{R}^n$ induced by  the  scalar product $g_0(x,y)=\omega(x,J_0y),$  $(x,0), (y,0) \in  \mathbb{R}^n \times \{0\} \simeq L.$ 
Elements of type $g_1, g_2$ and $g_3$ 
generate $Sp(V,\omega).$  See Folland \cite{Folland}. Note that 
in Habermann, Habermann \cite{HabermannHabermann}, a different convention for the Fourier transform is used.
Note that $L^2(\mathbb{R}^n)$ decomposes into the direct sum 
$L^2(\mathbb{R}^n)_+\oplus L^2(\mathbb{R}^n)_-$ of irreducible $\widetilde{G}$-modules of the even and of the odd functions in $L^2(\mathbb{R}^n).$  
For a proof that $\sigma$ is a representation,
see Folland \cite{Folland} or  Wallach 
\cite{Wallach} for instance. For a proof that $\sigma$ intertwines the 
Schr\"odinger representation of 
the Heisenberg group, see Wallach \cite{Wallach}, Habermann, Habermann 
\cite{HabermannHabermann} or  Folland \cite{Folland}. A proof that $L^2(L)_{\pm}$ are irreducible is contained  
in Folland \cite{Folland}.

\bigskip

Taking the derivative $\sigma_*$ at the unit element of $\widetilde{G}$  of the representation 
$\sigma$ restricted  to 
 smooth vectors in $L^2(L),$
we get the representation $\sigma_*: \widetilde{\mathfrak{g}} \to 
\mbox{End}\,(S) $ of the  Lie algebra of the metaplectic group
on the vector space $S=S(L)$ of Schwartz functions on $L.$ See Borel, Wallach 
\cite{BorelWallach} and
Folland \cite{Folland} where the smooth vectors are determined. Note that, we have 
$S \simeq S_+\oplus S_-$  similarly as in the previous decomposition.
For $n \times n$ real matrices $B=B^t,$ $C=C^t$ and $A,$  we have (see Folland \cite{Folland})
\begin{eqnarray*}
&&\sigma_*(X) =   \frac{1}{4\pi \imath} \sum_{i,j=1}^n B_{ij} 
\frac{\partial^2}{\partial x^i \partial x^j} 
\mbox{ for }  X=\left( \begin{array}{cc}
0 & B \\
0 & 0
\end{array} \right)\\
&&\sigma_*(Y) = - \pi \imath \sum_{i,j=1}^n C_{ij}  x^i x^j 
\mbox{ for }  Y=\left( \begin{array}{cc}
0 & 0 \\
C & 0
\end{array} \right)\\
&&\sigma_*(Z) =  - \sum_{i,j = 1}^{n} A_{ij}x^j\frac{\partial}{\partial x^i} - 
\frac{1}{2}\sum_{i=1}^n A_{ii}
\mbox{ for }  Z = \left(\begin{array}{cc}
A   &  0 \\
0 & -A^t
\end{array} \right). \\
\end{eqnarray*}	

\bigskip

It follows that 
$$\sigma_*(J_0) =  
\imath \sum_{i=1}^n\left(\frac{1}{4\pi}\frac{\partial^2}{\partial(x^i)^2} -  
\pi(x^i)^2\right).$$

\bigskip

\noindent{\bf Definition 1:} For any $m\in \mathbb{N}_0,$ we set
$h_m(x) = \frac{2^{1/4}}{\sqrt{m!}}(\frac{-1}{2\sqrt{\pi}})^me^{\pi x^2} 
\frac{d^m}{d x^m} (e^{-2\pi x^2}).$
For $n\in \mathbb{N}$ and $\alpha \in \mathbb{N}_0^n,$ we define the Hermite 
function $h_{\alpha}$ with 
index $\alpha=(\alpha_1,\ldots, \alpha_n)$ by $h_{\alpha}(x^1,\ldots,x^n) = 
h_{\alpha_1}(x^1)\ldots h_{\alpha_n}(x^n),$
$(x^1,\ldots, x^n) \in \mathbb{R}^n.$ 

\bigskip

\noindent{\bf Remark:} For Hermite functions, see Whittaker, Watson 
\cite{WhittakerWatson} and Folland \cite{Folland}. We use the convention of 
Folland 
\cite{Folland}. Especially, $h_0(x)=2^{1/4}e^{-\pi x^2}.$

\bigskip

Well known properties of Hermite functions make us able to derive that for any
$\alpha =(\alpha_1,\ldots$ $\ldots, \alpha_n) \in \mathbb{N}_0^n$
$$\sigma_*(J_0)h_{\alpha} = - \imath(|\alpha| + 
\frac{n}{2})h_{\alpha}$$
where 
$|\alpha|=\sum_{i=1}^n \alpha_i.$ Thus, the Hermite functions are eigenfunctions of $\sigma_*(J_0).$

\subsection{Weyl algebra and Symplectic spinor multiplication}

Let $\Bbbk$ be a field of characteristic zero.
For any $n\in \mathbb{N},$  the {\bf Weyl algebra} $W_n$ over $\Bbbk$
is the associative algebra generated by elements $1\in \Bbbk,$ $a_1,\ldots, 
a_n$ and 
$b_1,\ldots, b_n$ satisfying to the relations $1a_i=a_i1,$ $1b_i = b_i1,$ $a_i b_j 
- b_ja_i = -\delta_{ij} 1,$ $a_ia_j=a_ja_i,$ $b_ib_j=b_jb_i,$
$1\leq i,j \leq n.$ 
It is known that $W_n$ has a faithful representation on the space of 
polynomials 
$\Bbbk[q^1,\ldots, q^n]$ given by
$1 \mapsto 1$ (multiplication by $1$), $a_i \mapsto q^i$ and $b_i \mapsto 
\frac{\partial}{\partial q^i},$ where $q^i$ denotes the multiplication 
of a polynomial by $q^i$ and
$\frac{\partial}{\partial q^i}$ is the partial derivative in the $i$-th 
variable. See, e.g., Bj\"ork \cite{Bjork}.

\bigskip

Any associative algebra  $A$ over  field $\Bbbk$ can be equipped with the 
 commutator $$[\, ,\,]: A \times A \to A$$ defined by
$[x,y] = xy - yx, \mbox{ } x,y \in A,$ making it a Lie algebra. 
The {\bf Heisenberg Lie algebra} $H_n$ is the real vector space 
$\mathbb{R}^{2n+1}[q^1,\ldots, q^n, p_1,\ldots, p_n,t]$ with the Lie bracket $$[\, ,\,]: H_n \times H_n \to 
\{0\} \times \{0\} \times \mathbb{R} \subseteq H_n$$ prescribed on basis by 
$[\partial_t,\partial_{q^i}]=[\partial_t,\partial_{p_i}]=[\partial_{q^i}, 
\partial_{q^j}]
 = [\partial_{p_i},\partial_{p_j}] = 0$ and 
$[\partial_{q^i},\partial_{p_j}] = -\delta_{ij}\partial_t,$ $1\leq i, j\leq n.$ 
Note that   $[\, ,\,]$ is not the Lie bracket of vector fields in this case.
It is the 
Lie algebra of the Heisenberg group $H(n)$  and
it is isomorphic (as a Lie algebra) to  $$W_n(1) = \{\tau 1 +\sum_{i=1}^n 
(\alpha_i a_i + \beta_i b_i)|\, \tau, \alpha_i, \beta_i \in \mathbb{R}, \, 
i=1,\ldots, n\} 
\subseteq W_n$$  equipped with the commutator as the Lie algebra bracket.
An isomorphism can be given on a  basis by 
$\partial_t \mapsto 1,$ $\partial_{q^i}  \mapsto a_i,$  $\partial_{p_i} \mapsto 
b_i,$ $i=1,\ldots, n.$

\bigskip

For a symplectic vector space $(V,\omega)$ of dimension $2n$ over $\mathbb{R},$ 
let us choose a symplectic basis $(e_i)_{i=1}^{2n}$ and consider the tensor algebra  $$A=T(V^{\mathbb{C}}) = \mathbb{C} \oplus V^{\mathbb{C}} \oplus (V^{\mathbb{C}} \otimes 
V^{\mathbb{C}}) \oplus \cdots .$$ Let us set
$sCliff(V,\omega)= A / I,$
where $I$ is the two sided ideal  generated over $A$ by elements 
$v \otimes w - w \otimes v + \imath \omega(v,w),$ $v,w\in V^{\mathbb{C}}.$ 
The complex associative algebra $sCliff(V,\omega)$ is called the {\bf symplectic Clifford 
algebra}.
Let us consider the map  $1\mapsto 1,$ $e_{n+i} \mapsto - a_i$  and 
$e_{n+1-i} \mapsto \imath b_i,$ $i=1,\ldots, n,$ which extends to a 
homomorphism of associative algebras $sCliff(V,\omega)$ and $W_n$ for $\Bbbk = \mathbb{C}.$  It is not 
difficult to see that this map is an isomorphism onto $W_n.$ 
Summing up, $W_n$ and $sCliff(V,\omega)$ are isomorphic as associative algebras.
The Heisenberg Lie algebra  $H_n$ embeds homomorphically into 
$sCliff(V,\omega)$ (considered as a Lie 
algebra with respect to the commutator) via $\partial_t \mapsto 1,$ 
$\partial_{q^i} \mapsto -e_{n+i}$ and $\partial_{p_i} \mapsto \imath 
e_{n+1-i},$ 
$i=1,\ldots, n.$

\bigskip

\noindent{\bf Remark:} Note that there is an isomorphism of the Heisenberg Lie algebra $H_n$ with $\Bbbk_1[q^1,\ldots,$ $\ldots q^n, 
p_1,\ldots, p_n],$ the space of degree one polynomials in $q^i, p_i$ 
($i=1,\ldots, 
n$), equipped with the Poisson bracket $$\{f,g\}_P = \sum_{i,j=1}^n 
\left(\frac{\partial f}{\partial q^i} \frac{\partial g}{\partial p_i} - 
\frac{\partial f}{\partial p_i}\frac{\partial g}{\partial q^i}\right)$$
where $f,g \in \Bbbk_1[q^1, \ldots, q^n, p_1, \ldots, p_n].$

\bigskip

We come to the following important definition.

\bigskip

\noindent{\bf Definition 2:} Let $(e_i)_{i=1}^{2n}$ be a symplectic basis of $(V,\omega)$. 
For $i=1,\ldots, 
n$ and $f\in S,$ we set
$$e_i\cdot f =\imath x^i f \mbox{   \, \,  and  \, \,   } e_{i+n} \cdot f =\frac{\partial f}{
\partial x^{n-i+1}}$$ and extend it linearly to $V.$ The map $\cdot : V\times 
S \to S$ is called the {\bf symplectic spinor 
multiplication}.

\bigskip

\noindent{\bf Remark:} In the preceding definition, $f\in S(\mathbb{R}^n)$ and 
$x^i$ denotes the projection onto
the $i$-th coordinate in $\mathbb{R}^n.$ 
Note that the symplectic spinor multiplication depends on the choice of a 
symplectic basis. Because of its equivariant properties (see Habermann \cite{HabermannHabermann}, p. 13), one can use 
the multiplication on the level of bundles. In this case, we denote it by the dot 
as well. Note that the equivariance of the symplectic Clifford multiplication with 
respect to the Segal--Shale--Weil representation  makes the definitions of the 
symplectic spinor Dirac, the second symplectic spinor Dirac and the associated operator 
correct.



\section{Symplectic spinors in differential geometry}

Let us recall that a {\bf symplectic manifold} is a   manifold equipped with a 
closed non-degenerate exterior differential $2$-form $\omega.$ 

\bigskip

One of the big achievements of Bernhard Riemann in geometry is a
definition of the curvature (Kr\"ummungsma{\ss}) in an arbitrary dimension.
After publishing of his Habilitationsschrift,  Levi-Civita and Riemannian connections became fundamental objects for metric geometries.
Intrinsic notions and properties (such as straight lines, angle deficits,  parallelism etc.)  of many geometries known in that time  
can be defined and investigated by means of them. Using these connections,  one can find out quite easily, whether
the given manifold is locally isometric to the Euclidean space.

\bigskip

\noindent{\bf Definition 3:} Let $(M,\omega)$ be a symplectic manifold. An affine 
connection $\nabla$ on $M$ is called symplectic
if $\nabla \omega = 0.$
Such a connection is called a {\bf Fedosov connection} if it is torsion-free.

\bigskip

For  symplectic connections, see, e.g.,  Libermann \cite{Libermann}, Tondeur 
\cite{Tondeur}, Vaisman \cite{Vaisman} and Gelfand, Retakh, Shubin \cite{GRS}.
In contrast to Riemannian geometry, we have the following theorem which goes back to Tondeur \cite{Tondeur}. See  Vaisman
\cite{Vaisman} for a proof.

\bigskip

\noindent{\bf Theorem 4}: The space of Fedosov connections on a symplectic manifold $(M,\omega)$ is isomorphic to 
an affine space modeled on the infinite dimensional vector space  $\Gamma(S^3 TM),$ where $S^3 TM$ denotes the 
third symmetric product of the tangent bundle of $M.$

\bigskip

\noindent{\bf Remark:} Note that due to a theorem of Darboux (see McDuff, Salamon \cite{McDuffSalamon}), all symplectic manifolds of 
equal dimension are locally equivalent. In particular, symplectic connections 
cannot serve for  distinguishing of symplectic manifolds in the local sense.
From the eighties of the last century, symplectic 
connections gained an important role in mathematical physics. They became 
crucial for   quantization procedures. See Fedosov 
\cite{Fedosov} and Waldmann \cite{Waldmann}. 

\bigskip

Let $(M^{2n},\omega)$ be a symplectic manifold and $\nabla$ be a Fedosov connection. 
The {\bf curvature tensor field} of $\nabla$ is defined by $$R(X,Y)Z= 
\nabla_X \nabla_Y Z -\nabla_Y\nabla_XZ - \nabla_{[X,Y]}Z$$ where $X,Y,Z \in 
\mathfrak{X}(M).$ A {\bf local symplectic frame} $(U,(e_i)_{i=1}^{2n})$ of 
 $(M,\omega)$ is an open subset $U$ in $M$ and a  sequence of vector fields
$e_i$ on $U$ such that $((e_i)_m)_{i=1}^{2n}$ is a symplectic basis of
$(T_mM,\omega_m)$ for each $m \in U.$

\bigskip

Let $(U,(e_i)_{i=1}^{2n})$ be a local symplectic
frame. For $X=\sum_{i=1}^{2n}X^ie_i, Y=\sum_{i=1}^{2n}Y^ie_i, Z = \sum_{i=1}^{2n}Z^ie_i, 
V=\sum_{i=1}^{2n}V^ie_i \in \mathfrak{X}(M),$ $X^i, Y^i,Z^i, V^i \in 
\mathcal{C}^{\infty}(U),$  and 
$i,j,k,l = 1, \ldots, 2n,$ we set
\begin{eqnarray*}
R_{ijkl} &=& \omega(R(e_k,e_l)e_j,e_i)\\
\sigma(X,Y) &=& \mbox{Tr}(V \mapsto R(V,X)Y), \, V \in \mathfrak{X}(M)\\
\sigma_{ij} &=&\sigma(e_i,e_j)\\
\sigma_{ijkl}&=& \frac{1}{2(n+1)}(\omega_{il}\sigma_{jk} 
-\omega_{ij}\sigma_{lk} + \omega_{jl}\sigma_{ik} 
-\omega_{jl} \sigma_{ik} + 2\sigma_{ij}\omega_{kl})\\
\widetilde{\sigma}(X,Y,Z,V) &=& 
\sum_{i,j,k,l=1}^{2n}\sigma_{ijkl}X^iY^jZ^kV^l \\
W&=& R - \widetilde{\sigma}
\end{eqnarray*}
where at the last row, $R$ represents the $(4,0)$-type tensor field 
$\sum_{i,j,k,l=1}^{2n}R_{ijkl}\epsilon^i\otimes \epsilon^j \otimes \epsilon^k \otimes \epsilon^l$ and
$(\epsilon^i)_{i=1}^{2n}$ is the frame dual to $(e_i)_{i=1}^{2n}.$

\bigskip

\noindent{\bf Definition 4:} We call $W$ the {\bf symplectic Weyl 
curvature}. The  $(2,0)$-type   tensor field $\sigma$
is called the {\bf symplectic Ricci curvature}. A symplectic manifold with a Fedosov
connection is called of Ricci-type if $W=0$ and it is called Ricci-flat if $\sigma = 0.$

\bigskip

Let $(M,\omega)$ be a symplectic manifold. We set
$$Q=\{f  \mbox{ is a symplectic basis of } 
(T_mM,\omega_m)|\, m\in M\}$$ and call it the {\bf symplectic rep\`ere bundle}. For any $f=(e_1,\ldots, e_{2n}) \in 
Q,$ we denote by 
$\pi_Q(f)$ the unique point $m\in M$ such that each vector in $f$ belongs to 
$T_mM.$ 
The topology on $Q$ is the coarsest  one for which $\pi_Q$ is continuous.
It can be seen  that $\pi_Q:Q \to M$ is a  
principal $Sp(2n,\mathbb{R})$-bundle.   

\bigskip 

\noindent{\bf Definition 5:}
A pair $(P,\Lambda)$ is called a {\bf metaplectic structure} if $\pi_P:P\to M$ is 
a  principal $Mp(2n,\mathbb{R})$-bundle over $M$ and $\Lambda:P \to Q$ is a 
principal bundle homomorphism such that
the following diagram commutes. The horizontal arrows denote the actions of 
$\widetilde{G}$ and $G,$ respectively.
$$\begin{xy}\xymatrix{
P \times \widetilde{G} \ar[dd]^{\Lambda\times \lambda} \ar[r]&   
P \ar[dd]^{\Lambda} \ar[dr]^{\pi_P} &\\
                                                            & &M\\
Q \times G \ar[r]   & Q \ar[ur]_{\pi_Q} }\end{xy}$$

\bigskip

 A {\bf compatible positive almost complex structure} $J$ on  a symplectic manifold 
$(M, \omega)$ is any endomorphism
$J: TM \to TM$ such that $J^2=-1_{TM}$ and such that 
$g(X,Y) = \omega(X, JY),$ $X,Y \in \mathfrak{X}(M),$ is a Riemannian metric. In 
particular, $g$ is a symmetric tensor field. Note that $J$ is an isometry 
and a symplectomorphism as well. A compatible positive
almost complex structure always exists on a symplectic manifold $(M,\omega).$
See, e.g., McDuff, Salamon \cite{McDuffSalamon}, pp. 63 and 70, for a proof.

\bigskip 

\noindent{\bf Remark:} Note that a K\"ahler manifold can be defined as  a symplectic manifold equipped with a Fedosov connection $\nabla$ and a
compatible positive almost complex structure $J$ such that $\nabla J = 0,$ i.e., $J$ is $\nabla$-flat. Especially, any K\"ahler manifold is symplectic. 
 The first example of a compact symplectic manifold which does not admit any K\"ahler structure was given by Thurston \cite{Thurston}.
He was inspired  by a review note of Libermann  \cite{Libermannreview} who comments a mistake in an article of Guggenheimer \cite{Guggenheimer}.
See also the review \cite{Hodge} of the Guggenheimer's article by Hodge.

\bigskip

In the following theorem, a condition  for the 
existence of a metaplectic structure is given.

\bigskip

\noindent{\bf Theorem 5:} Let $(M,\omega)$ be a symplectic manifold and $J$ be a 
compatible positive almost complex structure. Then $(M,\omega)$ possesses a 
metaplectic structure if and only if 
the second Stiefel-Whitney class  $w_2(TM)$ of $TM$ vanishes 
if and only if the first Chern class $c_1(TM) \in H^2(M,\mathbb{Z})$ of 
$(TM,J)$ is even.

\noindent{\it Proof.} See Kostant \cite{Kostant} and Forger, Hess \cite{ForgerHess}, p. 270. \hfill $\Box$

\bigskip

\noindent{\bf Remark:} An element $a\in H^2(M,\mathbb{Z})$ is called even if 
there is an element $b\in H^2(M,\mathbb{Z})$ such that $a=2b.$ By
a Chern class of $(TM,J),$ we mean the Chern class of the complexification
$TM^{\mathbb{C}}$ defined with the help of the compatible positive almost complex structure $J$. See 
Milnor, Stasheff \cite{MilnorStasheff}.

\subsection{Habermann's symplectic Dirac and associated second order operator}

We introduce the  symplectic Dirac operators  and the  associated 
second order operator of K. Habermann. Note that there exists a complex version
of the metaplectic structure (so-called $Mp^c$-structure), and also of the mentioned
operators.
See Robinson, Rawnsley \cite{RobinsonRawnsley} and Cahen,  Gutt, La Fuente Gravy 
and Rawnsley \cite{CGGR}. Let us notice that
$Mp^c$ structures exist globally on any symplectic 
manifold (see \cite{RobinsonRawnsley}). Generalizations of many results of Habermann, Habermann in \cite{HabermannHabermann}
to the $Mp^c$-case are  straightforward (see \cite{CGGR}).

\bigskip

\noindent{\bf Definition 6:} Let $(M^{2n},\omega)$ be a symplectic manifold admitting a 
metaplectic structure $(P,\Lambda).$ 
The associated bundle $\mathcal{S} = P \times_{\sigma} S$ is called  
the {\bf symplectic spinor} or the {\bf Kostant's bundle}.  Its smooth sections are 
called {\bf symplectic spinor fields}.


\bigskip 

After introducing the Kostant's bundle, we can set up  definitions of the differential operators.

\bigskip

\noindent{\bf Definition 7:} Let $\nabla$ be a symplectic connection on a symplectic 
manifold $(M,\omega)$ admitting a metaplectic structure $(P,\Lambda)$. Consider the principal 
connection $TQ \to \mathfrak{sp}(2n,\mathbb{R})$ induced
by $\nabla$ and its lift $Z: TP \to \widetilde{\mathfrak{g}}$ to the 
metaplectic structure.  The associated
covariant derivative $\nabla^S: \Gamma(\mathcal{S}) \to \Gamma(\mathcal{S} \otimes T^*M)$ on 
symplectic spinor fields is called the {\bf symplectic spinor covariant derivative}.  
Let $(U,(e_i)_{i=1}^{2n})$ be a local symplectic frame. The  operator 
$D: \Gamma(\mathcal{S}) \to \Gamma(\mathcal{S})$ 
defined for any $\phi \in \Gamma(\mathcal{S})$ by 
$$D\phi =  \sum_{i,j=1}^{2n} \omega^{ij}e_i \cdot \nabla^S_{e_j} \phi$$ 
is called the {\bf (Habermann's)  symplectic spinor Dirac operator}.

\bigskip

Let $J$ be a compatible positive almost complex structure on a symplectic manifold $(M,\omega).$
A {\bf local unitary frame} is a local symplectic frame 
which is 
orthogonal with respect
to the associated Riemann  tensor $g(X,Y) = \omega(X, JY),$ $X,Y \in 
\mathfrak{X}(M).$

\bigskip

\noindent{\bf Definition 8:}
Let $J$ be a compatible positive almost complex  structure on a symplectic manifold which admits a metaplectic structure and
$(U,(e_i)_{i=1}^{2n})$ be a local unitary frame.
The operator $\widetilde{D}:\Gamma(\mathcal{S}) 
\to \Gamma(\mathcal{S})$ defined for any $\phi \in \Gamma(\mathcal{S})$ by  $$\widetilde{D}\phi = \sum_{i = 1}^{2n}  
(Je_i)\cdot  \nabla^S_{e_i}\phi$$ 
is called the {\bf second symplectic spinor Dirac operator}.
The operator $\mathfrak{P} = \imath[\widetilde{D}, D]$ is called the 
{\bf associated second order operator}. 

\bigskip

\noindent{\bf Remark:}
The associated second order operator $\mathfrak{P}$ is elliptic in the sense that its principal symbol
$\sigma(\mathfrak{P},\xi): \mathcal{S} \to \mathcal{S}$ is a bundle isomorphism for any non-zero 
cotangent vector $\xi \in T^*M.$ See Habermann, Habermann 
\cite{HabermannHabermann}, p. 68.

\bigskip

For symplectic spinor covariant derivative $\nabla^S$ and a chosen compatible positive almost complex structure,
one defines the formal adjoint $(\nabla^S)^*: \Gamma(\mathcal{S} \otimes T^*M) \to \Gamma(\mathcal{S})$ of $\nabla^S.$ See Habermann, 
Habermann  \cite{HabermannHabermann}.

\bigskip

\noindent{\bf Definition 9}: The {\bf Bochner-Laplace operator} on symplectic 
spinors $\Delta^S:\Gamma(\mathcal{S}) \to \Gamma(\mathcal{S})$ is the 
composition  $\Delta^S = (\nabla^S)^* \circ \nabla^S.$

\bigskip

\noindent{\bf Definition 10:} The curvature tensor field $R^S$ on  symplectic spinors  induced by a Fedosov connection $\nabla$ is defined by 
$$R^S(X,Y)\phi = \nabla^S_X\nabla^S_Y\phi-\nabla^S_Y\nabla^S_X\phi - 
\nabla_{[X,Y]}^S\phi$$
where $X,Y \in \mathfrak{X}(M),$ $\phi \in \Gamma(\mathcal{S})$ and $\nabla^S$ is the symplectic spinor derivative.

\bigskip

In the next theorem,   a relation of the associated second order operator 
$\mathfrak{P}$ to the Bochner-Laplace operator $\Delta^S$ on symplectic spinors is 
described.  It is derived by K. Habermann, and it is a parallel to the well  
known Weitzenb\"ock's and Lichnerowicz's formulas for the Laplace operator of the de Rham differentials (Hodge-Laplace) and the Laplace operator of a Levi-Civita connection (Bochner-Laplace); and 
for the square of the Dirac operator and  the Laplace operator of a Lichnerowicz 
connection on  
spinors  (Lichnerowicz-Laplace), respectively. See, e.g., Friedrich \cite{Friedrich} for the latter formula. We present a version 
of the Habermann's theorem for  K\"ahler manifold. See Habermann, Habermann 
\cite{HabermannHabermann} for more general versions.

\bigskip

\noindent{\bf Theorem 6:} Let $(M,\omega, J)$ be a K\"ahler manifold and $(U,(e_i)_{i=1}^{2n})$ be a local unitary frame.
 Then for any $\phi \in 
\Gamma(\mathcal{S})$
$$\mathfrak{P}\phi = \Delta^S \phi + \imath \sum_{i,j=1}^{2n} (J e_i) \cdot e_j \cdot
R^S(e_i, e_j)\phi.$$
 
\noindent{\it Proof.} See Habermann, Habermann \cite{HabermannHabermann}. \hfill $\Box$

\bigskip

For complex manifolds of complex dimension one\footnote{i.e., Riemann 
surfaces}, Habermann obtains the following  
consequence of the formula in Theorem 6.

\bigskip

\noindent{\bf Theorem 7:} If $M$ is a Riemann  surface of genus $g \geq 2,$ 
$\omega$ is a volume form on $M,$ and $(P,\Lambda)$ is a metaplectic 
structure, then the kernel of the associated second order operator is trivial.

\noindent{\it Proof.} See Habermann \cite{HabermannManuscripta}. \hfill $\Box$

\bigskip

\noindent{\bf Remark:} 
In \cite{HabermannCMP} and \cite{HabermannManuscripta},  Habermann proves that for  $T^2$  ($g=1$) and the 
trivial metaplectic structure, the null space for 
$\mathfrak{P}$ is isomorphic to the Schwartz space $S=S(\mathbb{R}).$
In the case of the (trivial) metaplectic structure on the sphere, the 
kernel of the associated second order operator is  rather complicated. See Habermann \cite{HabermannManuscripta} or Habermann, Habermann \cite{HabermannHabermann}.
 In the case of genus $g=1$ and non-trivial metaplectic structures, the kernel of 
$\mathfrak{P}$ is trivial as well. For it, see Habermann \cite{HabermannManuscripta}.

\bigskip

For further results on spectra and null-spaces of the introduced operators, 
see  Brasch, Habermann, Habermann \cite{BHH},   Cahen, La 
Fuente Gravy, Gutt, Rawnsley  \cite{CGGR} and Korman \cite{Korman}.
The key features used are the Weitzenb\"ock-type formula (Theorem 6)
and an orthogonal decomposition of the Kostant's  
bundle. To our knowledge,
this   decomposition was used firstly by Habermann in this context. It is   
derived from  a $\widetilde{K}$-isomorphism between 
$L^2(\mathbb{R}^n)$ and the 
Hilbert orthogonal sum $\bigoplus_{m=0}^{\infty} \mathfrak{H}_m$ of the spaces $$\mathfrak{H}_m = \bigoplus_{\alpha, |\alpha|\leq 
m}\mathbb{C} h_{\alpha}, \mbox{ } m \in \mathbb{N}_0.$$ Recall that $\widetilde{K}$ denotes the preimage in the metaplectic group of the unitary
group $U(n)$ by the covering $\lambda.$ (See Habermann, Habermann \cite{HabermannHabermann}, p. 18 for a description of the isomorphism.)

\subsection{Quantization by symplectic spinors}

For a symplectic manifold $(M,\omega)$ and a smooth function $f$ on $M,$ we denote
by $X_f$ the vector field $\omega$-dual to $df,$ i.e., such a vector field for which $$\omega(X_f, 
Y) = (df)Y$$ for any  $Y\in \mathfrak{X}(M).$ 
It is called the  Hamiltonian vector field of $f.$
A vector field is called symplectic if its flow preserves the symplectic form. 
Any Hamiltonian vector field is symplectic but not vice versa. For a study of these notions, we refer to the  monograph
McDuff, Salamon \cite{McDuffSalamon}.
Note that in this formalism, a Poisson bracket of two smooth functions $f,g$ on $M$ is defined  by
$$\{f,g\}_P = \omega(X_f,X_g).$$

\bigskip
Let $(M,\omega)$ be  a symplectic manifold admitting a metaplectic structure.
For a symplectic vector field $Y$, let  $L_Y$ denote the Lie derivative on the sections of the Kostant's bundle in 
direction $Y.$  See Habermann, Klein \cite{HabermannKlein} and Kol\'a\v{r}, 
Michor, Slov\'ak \cite{KolarMichorSlovak}.

\bigskip

\noindent{\bf Definition 11:} Let $(M,\omega)$ be  a symplectic manifold admitting a metaplectic structure. For a smooth function $f$
on $M$, we define a map $\mathfrak{q}(f): \Gamma(\mathcal{S}) \to 
\Gamma(\mathcal{S})$ by $$\mathfrak{q}(f)\phi = -\imath \hslash 
L_{X_f}\phi$$ for any $\phi\in \Gamma(\mathcal{S}).$ We call $\mathfrak{q}: f \mapsto \mathfrak{q}(f)$  the {\bf Habermann's map}.

\bigskip
 
Due to the properties of $L_X,$ it is clear that $\mathfrak{q}$ maps into the vector space endomorphisms of $\Gamma(\mathcal{S})$
$$\mathfrak{q}: \mathcal{C}^{\infty}(M) \to \mbox{End}(\Gamma(\mathcal{S})).$$   Moreover, Habermann derives the following theorem.

\bigskip 

\noindent{\bf Theorem 8:} Let $(M,\omega)$ be a symplectic manifold admitting a metaplectic struture.
Then for any $f,g \in \mathcal{C}^{\infty}(M),$ the Habermann's  map 
satisfies 
$$[\mathfrak{q}(f), \mathfrak{q}(g)] = \imath \hslash \, \mathfrak{q}(\{f,g\}_P).$$

\noindent{\it Proof.} See Habermann, Habermann \cite{HabermannHabermann}.\hfill  $\Box$

\bigskip

\noindent{\bf Remark:}
The Habermann's map $\mathfrak{q}$ satisfies 
the quantization condition (see Introduction) and thus, it gives an example of a 
non-deformed quantization. By this we mean that $\mathfrak{q}$ is a morphism of 
Poisson algebras $(\mathcal{C}^{\infty}(M), \, \{,\})$  and $(\mbox{End}\, (\Gamma(\mathcal{S})),\, [,])$ up to a multiple.
However  notice that usually, a quantization is 
demanded to be  a map on smooth functions $\mathcal{C}^{\infty}(M)$ 
defined on the phase space $M$ into the space of operators on the vector space $L^2(N)$ of $L^2$-functions  or $L^2$-sections of a  line bundle over $N$
where $N$ denotes the Riemannian manifold of the configuration space.  
 See Souriau 
\cite{Souriau} and Blau \cite{Blau} for conditions on quantization maps, their 
constructions and examples.

\section{Author's results in Symplectic spinor geometry}

We present  results achieved by the author in  differential geometry concerning
 symplectic spinors that we consider important and relevant.
We start with a chapter on representational theoretical, or if we wish equivariant, properties of exterior differential forms with values in symplectic spinors.

\subsection{Decomposition of tensor products and a Howe-type duality}

Let $\mathfrak{g}$ be the Lie algebra of symplectic group 
$Sp(2n,\mathbb{R})$, $\mathfrak{g}^{\mathbb{C}}$ the complexification of 
$\mathfrak{g},$ $\mathfrak{h}$ 
a Cartan subalgebra of $\mathfrak{g}^{\mathbb{C}},$ $\Delta^+$ a choice of 
positive roots, and $\{\varpi_i\}_{i=1}^n$ 
the 
set of fundamental weights with respect to these choices.
Let us denote the irreducible complex highest weight module with highest weight 
$\lambda\in \mathfrak{h}^*$ by $L(\lambda).$ For any $\lambda = \sum_{i=1}^n\lambda_i\varpi_i,$
we set   $L(\lambda_1, \ldots, \lambda_n) = L(\lambda).$ 
For $i=0,\ldots, 2n,$ we denote by $\sigma^i$ the tensor product representation
of the complexified symplectic Lie algebra $\mathfrak{g}^{\mathbb{C}}$ on $E^i=\bigwedge^i V^* \otimes S,$ i.e.,
$\sigma^i: \mathfrak{g}^{\mathbb{C}} \to \mbox{End}\, (E^i)$ and
$\sigma^i(X)(\alpha \otimes s) = \lambda_*^{\wedge i}(X)\alpha \otimes s+ \alpha \otimes \sigma_*(X) s$
for any $X \in \mathfrak{g}^{\mathbb{C}}, \alpha \in \bigwedge^i V^*$ and $s \in S,$ where 
$\lambda^{ \wedge i}_*$ denotes the action
of $\mathfrak{g}^{\mathbb{C}}$ on $\bigwedge^i V^*.$
We consider $E=\bigoplus_{i=0}^{2n} E^i$ equipped with  the direct sum representation $\sigma^{\bullet}(X)=(\sigma^0(X),\ldots, \sigma^{2n}(X)),$ 
$X\in \mathfrak{g}^{\mathbb{C}}.$ Let us notice that here, $\sigma_*$ denotes the complex linear extension
 of the representation $\sigma_*:\mathfrak{g} \to \mbox{End}(S)$ considered above.

\bigskip

\noindent{\bf Remark:} Note that there is a misprint in Kr\'ysl \cite{KryslJOLT2}. 
Namely,  the ``action'' of
$\mathfrak{g}$ on $E$ (denoted by $\mathbb{W}$ there) is prescribed  by $X(\alpha\otimes s) = \lambda_*^{\wedge 
i}(X)\alpha \otimes \sigma_*(X)s$ for $X\in \mathfrak{g},$
$\alpha \in \bigwedge^i V,$ $s\in S,$ and $i=0,\ldots, 2n.$ 
Actually, we meant the standard tensor product representation as given above, i.e.,
$X(\alpha\otimes s) = \lambda_*^{\wedge i}(X)\alpha \otimes s + \alpha \otimes 
\sigma_*(X)s$. However, the results in \cite{KryslJOLT2} are derived  for the correct action $\sigma^{\bullet}$
defined above.

\bigskip

\noindent{\bf Definition 12:}
Let us set
$\Xi = \{(i,j_i)|\, i = 0,\ldots, n,  j_i =  0,\ldots, i \} \cup 
\{(i,j_i)|\, i = n+1, \ldots, 2n,  j_i = 0, \ldots,  2n - i\},$  
 $\mbox{sgn}(+)=0,$  $\mbox{sgn}(-) = 1,$ and 
$$E^{ij}_{\pm}=L(\underbrace{\frac{1}{2},\cdots,\frac{1}{2}}_{j},
\underbrace{-\frac{1}{2},\cdots,
-\frac{1}{2}}_{n-j-1},-1+\frac{1}{2}(-1)^{i+j+\mbox{sgn}(\pm)})$$ for 
$i=0,\ldots,n-1,$   $j=0,\ldots, i$ and $i=n,$   $j=0,\ldots, n-1.$ For $i=j=n,$
we set $E_+^{nn}=L(\frac{1}{2}, \cdots, \frac{1}{2})$ and 
$E_-^{nn}=L(\frac{1}{2}, \cdots, \frac{1}{2}, -\frac{5}{2}).$ 
For $i=n+1,\ldots,2n$ and $j=0,\ldots, 2n-i,$ we set 
$E^{ij}_{\pm}=E^{(2n-i)j}_{\pm}.$  
For any $(i,j) \in \mathbb{Z} \times \mathbb{Z} \setminus \Xi,$ we define $E^{ij}_{\pm}=0.$ 
Finally for any $(i,j)\in \mathbb{Z} \times \mathbb{Z}$, we set $E^{ij}=E^{ij}_+\oplus E^{ij}_-.$  
For $(i,j)\in \Xi$, the $\mathfrak{g}^{\mathbb{C}}$-modules $E^{ij}_{\pm}$  are called
{\bf higher symplectic spinor modules} and their elements  {\bf higher symplectic spinors}.

\bigskip

\noindent{\bf Theorem 9:} The following decomposition into irreducible 
$\mathfrak{g}^{\mathbb{C}}$-modules  holds
$$\bigwedge^{\bullet} V^* \otimes S_{\pm} = \bigoplus_{(i,j)\in \Xi} E^{ij}_{\pm}.$$

\noindent{\it Proof.} Kr\'ysl \cite{KryslJOLT2}. \hfill $\Box$

\bigskip

\bigskip

\noindent{\bf Remark:} The decomposition holds also on the level of minimal and hyperfunction 
 globalizations   since the corresponding  
globalization functors are adjoint functors to the 
Harish-Chandra forgetful functor. See Vogan \cite{Vogan} and Casselmann \cite{Casselmann}. It holds also for smooth Fr\'echet globalization $\widetilde{G} \to \mbox{Aut}(S).$
By abuse of notation, we shall denote the tensor product representation of $\widetilde{G}$ on $E$ by $\sigma^{\bullet}$ as well.
 The above decomposition holds also  when $V^*$ is replaced by $V$  
since the symplectic form gives an isomorphism of the appropriate 
representations of $\mathfrak{g}^{\mathbb{C}}.$

\bigskip

\noindent{\bf Definition 13:} For $i=0,\ldots, 2n,$ we denote the uniquely determined equivariant  projections of 
$\bigwedge^i V \otimes S^{\pm} \to E^{ij}_{\pm} \subseteq \bigwedge^i V 
\otimes S_{\pm}$ by $p_{\pm}^{ij}$ and the projections $p^{ij}_+ + p^{ij}_-$ onto $E^{ij}$ by 
$p^{ij},$ $(i,j)\in \mathbb{Z} \times \mathbb{Z}.$

\bigskip

 Let us recall a definition of the  simple Lie superalgebra $\mathfrak{osp}(1|2).$ It is 
generated by elements $e^+, e^-, h, f^+, f^-$   satisfying the following 
relations
\begin{equation*}
\begin{split}
[h,e^{\pm}] &= \pm e^{\pm}  \\
[h,f^{\pm}] &= \pm \frac{1}{2}f^{\pm} \\
[e^{\pm},f^{\mp}] &= -f^{\pm} 
\end{split}
\quad \quad \quad
\begin{split}
[e^+,e^-] &= 2h \\
\{f^{+},f^-\} &= \frac{1}{2}h\\
\{f^{\pm},f^{\pm}\} &= \pm \frac{1}{2}e^{\pm}
\end{split}
\end{equation*}
where $\{,\}$ denotes the anticommutator, i.e., $\{a,b\}=ab+ba,$ $a,b \in \mathfrak{osp}(1|2).$
\bigskip

We give a $\mathbb{Z}_2$-grading to the vector space $E=\bigwedge^{\bullet}V \otimes S$ by setting
$E_0 = \bigoplus_{i=0}^n\bigwedge^{2i} V \otimes S,$  $E_1 = \bigoplus_{i=1}^n\bigwedge^{2i-1} V 
\otimes S$ and $E=E_0\oplus E_1.$ Further, we choose a symplectic basis $(e_i)_{i=1}^{2n}$ of $(V,\omega)$ 
and denote its dual basis by $(\epsilon^i)_{i=1}^{2n} \subseteq V^*.$
The Lie superalgebra $\mathfrak{osp}(1|2)$ has a
representation $\rho: \mathfrak{osp}(1|2) \to \mbox{End}(E)$ on the superspace $E$ given by
\begin{equation*}
\begin{split}
\rho(f^{+})(\alpha \otimes s)= \frac{\imath}{2}\sum_{i=1}^{2n} 
\epsilon^i \wedge \alpha \otimes e_i \cdot s
\end{split}
\quad \mbox{and}   \quad 
\begin{split}
\rho(f^-)(\alpha \otimes s) = 
\frac{1}{2}\sum_{i=1}^{2n}\omega^{ij}\iota_{e_i}\alpha 
\otimes e_j \cdot s
\end{split}
\end{equation*}
where $\alpha \in \bigwedge^{\bullet} V^*,$ $s\in S,$ and $\iota_v$ denotes the contraction by the vector $v.$
Consequently, elements $ e^+, e^-$ and $h$ act by
\begin{equation*}
\begin{split}
\rho(e^{\pm})=\pm 2\{\rho(f^{\pm}),\rho(f^{\pm})\}
\end{split} 
\quad   \mbox{and} \quad 
\begin{split}
\rho(h)  = \frac{1}{2}[\rho(e^+),\rho(e^-)]
\end{split}
\end{equation*}
where $\{\, , \,\}$  and $[\, , \,]$ denote the anticommutator and the commutator on the associative algebra $\mbox{End}(E),$ respectively.

\bigskip

The following theorem is parallel to the Schur and Weyl dualities for tensor 
representations of $GL(n,\mathbb{C})$  and $SO(n,\mathbb{C}),$ respectively. See Howe 
\cite{HoweInvariant} where they are treated.

\bigskip

\noindent{\bf Theorem 10:}
The following $\mathfrak{g}^{\mathbb{C}} \times \mathfrak{osp}(1|2)$-module isomorphism holds
$$\bigwedge^{\bullet} V^* \otimes S \simeq \bigoplus_{i=0}^n \left( E^{ii}_+ 
\otimes F_i \right) \oplus \bigoplus_{i=0}^n \left( E^{ii}_-\otimes F_i\right) $$ where
$F_i=\mathbb{C}^{2n-2i+1}$  and  
$\rho_i:\mathfrak{osp}(1|2) \to \mbox{End}(F_i)$  is
 given on a basis $(b_j)_{j=i}^{2n-i}$ of $F_i$ by prescriptions  
\begin{equation*}
\begin{split}
&\rho_i(f^+)(b_j) =  A(n,i+1,j)b_{j+1}\\ 
&\rho_i(h)=2\{\rho_i(f^+),\rho_i(f^-)\} \, \, \, \, \, \, \, \, \mbox{ and}
\end{split}
\quad \quad \quad 
\begin{split}
&\rho_i(f^-)(b_j) = b_{j-1}\\
&\rho_i(e^{\pm}) = \pm 2\{\rho_i(f^{\pm}), \rho_i(f^{\pm})\}
\end{split}
\end{equation*}
 where $i=0,\ldots, n$ and $A(n,i,j) = \frac{(-1)^{i-j}+1}{16}(j-i) + 
\frac{(-1)^{i-j+1}+1}{16}(i+j-2n-1).$

\noindent{\it Proof.} See Kr\'ysl \cite{KryslJOLT2}. \hfill $\Box$

\bigskip

\noindent{\bf Remark:} In the preceding definition, if an index exceeds its allowed range, the object is considered to be zero.
Thus, e.g., $b_{2n-i+1}$ or $b_{i-2}$ are zero vectors.

\bigskip

\noindent{\bf Theorem 11:} 
For $i=0,\ldots, n,$ representations $F_i$ are irreducible.

\noindent{\it Proof.} See Kr\'ysl \cite{KryslJOLT2}. \hfill $\Box$

\bigskip

\noindent{\bf Remark:} Representations $\rho_i$ in Theorem 10 depend on the choice of a basis, but not their equivalence class. 
As follows from Theorem 11, the multiplicity of $E^{ii}_{\pm}$ in the $\mathfrak{g}^{\mathbb{C}}$-module
$E$ is $2n-2i+1$ for $i=0,\ldots, n.$

\subsection{Differential geometry of higher symplectic spinors}

For any symplectic manifold $(M,\omega)$ admitting a metaplectic structure $(P,\Lambda)$,
the decomposition from Theorem 9 can be lifted to the associated bundle
$\mathcal{E} = P \times_{\sigma^{\bullet}} E.$

\bigskip

\noindent{\bf Remark:} Since $S$ is a smooth globalization, we may consider $E$ as a representation of the metaplectic group as well.

\bigskip

\noindent{\bf Definition 14:} Let $(M,\omega)$ be a symplectic manifold admitting a 
metaplectic structure $(P,\Lambda)$.
 For any $(i,j) \in \mathbb{Z} \times \mathbb{Z},$ we set
$\mathcal{E}^{ij}= P \times_{\widetilde{G}} E^{ij}$ and call it the 
{\bf higher symplectic  spinor bundle} and elements of its section spaces the
{\bf higher symplectic spinor fields} if $(i,j) \in \Xi.$

\bigskip

We keep denoting the lifts of the projections $\bigwedge^i V^* \otimes S \to E^{ij}$ to 
$\Gamma(\mathcal{E}^i)\to \Gamma(\mathcal{E}^{ij})$ by $p^{ij},$ where  $\mathcal{E}^i = P\times_{\sigma^i} E^i.$


\subsubsection{Curvature, higher curvature and symplectic twistor complexes}

For a Fedosov connection $\nabla$ on a symplectic manifold $(M,\omega)$
admitting a metaplectic structure, we  consider the exterior covariant derivative $d^{\nabla^S}$ for the induced symplectic spinor 
derivative $\nabla^S.$ See, e.g., Kol\'a\v{r}, Michor, Slov\'ak \cite{KolarMichorSlovak} for a general construction of such derivatives.

\bigskip

\noindent{\bf Theorem 12:} Let $(M,\omega)$ be a symplectic manifold admitting a 
metaplectic structure and $\nabla$ be a Fedosov connection.
Then for any $(i,j) \in \mathbb{Z} \times \mathbb{Z},$ the restriction of the 
exterior symplectic spinor derivative 
satisfies 
$$d^{\nabla^S}:\Gamma(\mathcal{E}^{ij})\to \Gamma(\mathcal{E}^{i+1, j-1})\oplus 
\Gamma(\mathcal{E}^{i+1,j})\oplus \Gamma(\mathcal{E}^{i+1,j+1}).$$

\noindent{\it Proof.} See Kr\'ysl \cite{KryslArchivumSSVF}. \hfill $\Box$

\bigskip

\noindent{\bf Remark:} In particular, sections of each higher symplectic spinor 
bundle
are mapped into sections of at most three higher symplectic spinor bundles.
Note that in the case of orthogonal spinors in pseudo-Riemannian geometry, the 
target space structure of the exterior covariant derivative is similar. See 
Slupinski 
\cite{Slupinski}.

\bigskip
 

Let $(e_i)_{i=1}^{2n}$ be a local symplectic frame on $(M,\omega)$ and 
$(\epsilon^i)_{i=1}^{2n}$ be its dual symplectic coframe. Recall that above,
 we defined the symplectic Ricci and symplectic Weyl curvature tensor fields.
Let us denote by $\sigma^S$ the  endomorphism of the 
symplectic spinor bundle defined for any $\phi \in \mathcal{S}$ by
$$\sigma^S \phi 
= \frac{\imath}{2} 
\sum_{i,j,k,l=1}^{2n}{\sigma^{ij}}_{kl}\epsilon^k \wedge 
\epsilon^l \otimes e_i \cdot e_j \cdot \phi.$$

Similarly we set
$$W^S \phi  = \frac{\imath}{2} 
\sum_{i,j,k,l=1}^{2n}{W^{ij}}_{kl}\epsilon^k \wedge 
\epsilon^l \otimes e_i \cdot e_j \cdot \phi.$$

\bigskip

Recall that $$\bigwedge^2 T^*M \otimes \mathcal{S} = \mathcal{E}^{20}\oplus 
\mathcal{E}^{21} \oplus \mathcal{E}^{22}$$ according to Theorem 9.

\bigskip

In the next theorem,  components of $R^S$ in $\mathcal{E}^{20},$ 
$\mathcal{E}^{21}$ and $\mathcal{E}^{22}$ are found. 
We notice that
\begin{itemize}
\item[1)] we use the summation convention, i.e., if two indices occur which 
are labeled by the same letter, we sum over it without denoting the sum 
explicitly and
\item[2)] instead of $e_i\cdot e_j\cdot $, we write $e_{ij}\cdot$ and 
similarly for a higher number of indices.
\end{itemize}

\bigskip

\noindent{\bf Theorem 13:} Let $n>1$, $(M^{2n},\omega)$ be a symplectic manifold 
admitting a metaplectic structure and $\nabla$
be a Fedosov connection.
Then for any $\phi \in \Gamma(\mathcal{S}),$ $\sigma^S\phi \in \Gamma(\mathcal{E}^{20}\oplus \mathcal{E}^{21})$ and 
$W^S\phi \in \Gamma(\mathcal{E}^{21} \oplus \mathcal{E}^{22}).$ Moreover, we have the following 
projection formulas
\begin{eqnarray*}
p^{20} R^S\phi &=& \frac{\imath}{2n} \sigma^{ij}\omega_{kl}\epsilon^k \wedge 
\epsilon^l 
\otimes e_{ij}\cdot \phi\\
p^{21} R^S\phi &=& \frac{\imath}{n+1}\sigma^{ij}\epsilon^k\wedge 
\epsilon^l\otimes (\omega_{il} 
e_{kj} \cdot - \frac{1}{2n}\omega_{kl}e_{ij} \cdot )\phi -\\
&&\frac{\imath}{1-n} {W^{ijk}}_l \epsilon^m \wedge \epsilon^l \otimes e_{mkij} 
\cdot 
\phi \\
p^{22} R^S\phi &=& \frac{\imath}{2}{W^{ij}}_{kl} \epsilon^k \wedge \epsilon^l 
\otimes e_{ij} \cdot
 \phi + \frac{\imath}{1-n}{W^{ijk}}_l \epsilon^m \wedge \epsilon^l\otimes 
e_{mkij} \cdot \phi.
\end{eqnarray*}

\noindent{\it Proof.} See Kr\'ysl \cite{KryslJGP1}. \hfill $\Box$

\bigskip

\noindent{\bf Remark:} Note that for $n=1,$ $E^{21} = E^{22} = 0,$ so that there is 
no Weyl component of the curvature tensor of a Fedosov connection in this dimension. 
The formula for $p^{20}$ holds also for $n=1.$

\bigskip

\noindent{\bf Definition 15:} For $(i,j), (i+1,k) \in \Xi,$  $a=0,\ldots,n-1 $ and 
$b=n,\ldots, 2n-1,$ let us set 
\begin{equation*}
\begin{split}
D^{ij}_{i+1,k} = p^{i+1,k}d^{\nabla^S}_{|\Gamma(\mathcal{E}^{ij})}:\Gamma(\mathcal{E}^{ij}) \to 
\Gamma(\mathcal{E}^{i+1,k})
\end{split}
,\quad
\begin{split}
T_a=D^{aa}_{a+1,a+1} 
\end{split}
\quad \mbox{ and } \quad
\begin{split} 
T_b = D^{b,2n-b}_{b+1,2n-b-1}. 
\end{split}
\end{equation*}
The operators $T_i,$ $i=0,\ldots,2n-1,$ are called  the {\bf symplectic twistor operators}.

\bigskip

Let $(V,\omega)$ be a symplectic vector space, $(e_i)_{i=1}^{2n}$ be a 
symplectic basis, $(\epsilon^i)_{i=1}^{2n}$ be a  basis of $V^*$ dual
 to $(e_i)_{i=1}^{2n},$ and
 $\sigma \in S^2 V^*$ be a bilinear form. For $\alpha \in \bigwedge^{\bullet}V^*$ and $s\in S,$ we set
 $$\Sigma^{\sigma}(\alpha \otimes s) = \sum_{i,j=1}^{2n}{\sigma^i}_j \epsilon^j \wedge \alpha 
\otimes e_i \cdot s$$ and
 $$\Theta^{\sigma}(\alpha \otimes s) = \sum_{i,j=1}^{2n}\alpha \otimes 
\sigma^{ij} e_i \cdot e_j \cdot s.$$
We keep denoting the corresponding tensors on symplectic spinor bundles  by the same symbols.  In this case, the
 the symplectic Ricci curvature tensor field plays the role of the tensor $\sigma.$

\bigskip

We use abbreviations 
$$E^{\pm} =\rho(e^{\pm}): E \to E \mbox{ and } F^{\pm} =\rho(f^{\pm}): E \to E.$$
Let  $(M,\omega)$ be a symplectic manifold which admits a metaplectic structure and $\nabla$ be a Fedosov connection of Ricci-type.
For a higher symplectic spinor field $\phi \in 
\Gamma(\mathcal{E}),$ we have (see Kr\'ysl \cite{KryslMonat}) the following formula
$$R^E \phi = 
\frac{1}{n+1}(E^+\Theta^{\sigma}  + 2 F^+ \Sigma^{\sigma})\phi.$$

\bigskip

\noindent{\bf Remark:} By the higher curvature,
 we understand  the curvature of $\nabla^S$ on higher symplectic spinors, i.e.,
  $R^E = d^{\nabla^S} \circ d^{\nabla^S}.$

\bigskip

The above formula is used for proving  the next theorem. 

\bigskip

\noindent{\bf Theorem 14:} Let $n>1,$ $(M^{2n},\omega)$ be
a symplectic manifold admitting a metaplectic structure and $\nabla$ be
a  Fedosov connection of Ricci-type. Then
$$0  \longrightarrow
  \Gamma(\mathcal{E}^{00}) \overset{T_0}{\longrightarrow} 
  \Gamma(\mathcal{E}^{11}) \overset{T_{1}}{\longrightarrow}  
  \cdots  \overset{T_{l-1}}{\longrightarrow}
  \Gamma(\mathcal{E}^{nn}) \longrightarrow 0 \mbox{   and}  
$$ 
$$0  \longrightarrow
  \Gamma(\mathcal{E}^{nn}) \overset{T_n}{\longrightarrow} 
  \Gamma(\mathcal{E}^{n+1,n+1}) \overset{T_{n+1}}{\longrightarrow}  
  \cdots  \overset{T_{2n-1}}{\longrightarrow}
  \Gamma(\mathcal{E}^{2n,2n}) \longrightarrow 0  
$$ 
are complexes.

\noindent{\it Proof.} See Kr\'ysl \cite{KryslMonat}. \hfill $\Box$

\bigskip

We call the complexes from Theorem 14 the {\bf symplectic twistor complexes}.

\bigskip

\noindent{\bf Theorem 15:} Let $n>1,$ $(M^{2n},\omega)$ be a symplectic manifold admitting a 
metaplectic structure and $\nabla$ be a Fedosov connection of 
Ricci-type.  Then
$$ 0\longrightarrow \Gamma(\mathcal{E}^{00}) 
\overset{T_0}{\longrightarrow}\cdots \overset{T_{n-2}}{\longrightarrow} 
\Gamma(\mathcal{E}^{n-1,n-1})\overset{T_nT_{n-1}}{\longrightarrow}
\Gamma(\mathcal{E}^{n+1,n+1}) 
\overset{T_{n+1}}{\longrightarrow}
\cdots \overset{T_{2n-1}}{\longrightarrow} 
\Gamma(\mathcal{E}^{2n,2n}) \longrightarrow 0$$
is a complex.

\noindent{\it Proof.} See Kr\'ysl \cite{KryslMonat}. \hfill $\Box$

\bigskip

\noindent{\bf Definition 16:} Let $(\mathcal{F}^i \to M)_{i \in \mathbb{Z}}$ be a sequence of vector bundles over a smooth manifold $M,$
 $D^{\bullet}=(\Gamma(\mathcal{F}^i), D_i:\Gamma(\mathcal{F}^i) \to \Gamma(\mathcal{F}^{i+1}))_{i\in \mathbb{Z}}$ be a  complex of pseudodifferential operators and
for each $\xi \in T^*M,$ let $\sigma(D)(\xi)^{\bullet}=(\mathcal{F}_i, \sigma(D_i, \xi): \mathcal{F}^i \to \mathcal{F}^{i+1})_{i\in \mathbb{Z}}$ be the complex of symbols evaluated in $\xi$ which is
associated to the complex $D^{\bullet}.$ We call $D^{\bullet}$ {\bf elliptic} if
$\sigma(D)(\xi)^{\bullet}$ is an exact sequence in the category of vector bundles for any $\xi \in T^*M\setminus \{0\}.$

\bigskip

\noindent{\bf Remark}: Note that in  homological algebra, the above complexes are usually called cochain complexes.

\bigskip

\noindent{\bf Theorem 16:} Let $n>1,$  $(M^{2n},\omega)$ be a symplectic manifold admitting a metaplectic structure and
$\nabla$ be a  Fedosov connection of Ricci-type. 
Then the  complexes 
$$0  \longrightarrow
  \Gamma(\mathcal{E}^{0}) \overset{T_0}{\longrightarrow} 
  \Gamma(\mathcal{E}^{1}) \overset{T_{1}}{\longrightarrow}  
  \cdots  \overset{T_{n-2}}{\longrightarrow}  \Gamma(\mathcal{E}^{n-1}) \mbox{ 
  and}  
$$ 
$$
  \Gamma(\mathcal{E}^{n}) \overset{T_n}{\longrightarrow} 
  \Gamma(\mathcal{E}^{n+1}) \overset{T_{n+1}}{\longrightarrow}  
  \cdots  \overset{T_{2n-1}}{\longrightarrow}
  \Gamma(\mathcal{E}^{2n}) \longrightarrow 0  
$$ 
 are elliptic.
 
\noindent{\it Proof.} See Kr\'ysl \cite{KryslArchivumEC}. \hfill $\Box$

\subsubsection{Symplectic spinor Dirac, twistor and Rarita--Schwinger operators}

\bigskip

\noindent{\bf Definition 17:} Let $(M,\omega)$ be a symplectic manifold admitting  a metaplectic structure and
$\nabla$ be a Fedosov connection. The operators
\begin{equation*}
\begin{split}
\mathfrak{D} =  F^- \circ D^{00}_{10}:\Gamma(\mathcal{S}) \to 
\Gamma(\mathcal{S}) 
\end{split}
\quad \mbox{ \, \, and }
\begin{split}
\mathfrak{R} =  F^- \circ D^{11}_{21}: \Gamma(\mathcal{E}^{11}) \to 
\Gamma(\mathcal{E}^{11})
\end{split}
\end{equation*}
are called the {\bf symplectic spinor Dirac} and the {\bf symplectic spinor Rarita--Schwinger operator}, respectively.

\bigskip

\noindent{\bf Remark:}
$\mathfrak{D}$ is the $1/2$ multiple of the Habermann's symplectic spinor Dirac operator. 
 
\bigskip

Let us denote the set of  eigenvectors of a vector space endomorphism 
$G: W \to W$  by $\mbox{eigen}(G)$ and the 
set of its eigenvalues by ${\it spec}(G).$  Recall that by an eigenvalue, we mean simply a 
complex number $\mu,$ for which there is a nonzero $w\in W,$ such that 
$Gw=\mu w.$  (We do not investigate spectra from  the functional
analysis  point of view.)

\bigskip

\noindent{\bf Definition 18:}
A {\bf symplectic Killing spinor field} is any not everywhere zero section 
$\phi\in \Gamma(\mathcal{S})$ for which there exists $\mu\in \mathbb{C}$ such 
that
$$\nabla^{S}_X\phi=\mu X \cdot \phi$$ for each $X \in \mathfrak{X}(M).$ (The 
dot denotes the symplectic Clifford multiplication.) The set 
of symplectic  Killing spinor fields is denoted by $\mbox{kill}.$ Number 
$\mu$ from the above equation is called the {\bf symplectic Killing spinor 
number} 
and its set is denoted by $\it{kill}.$

\bigskip

\noindent{\bf Remark:}  The equation for a symplectic Killing spinor field can be 
written also 
as  $$\nabla^S \phi  = - 2 \mu \imath F^+ \phi.$$

\bigskip

\noindent{\bf Remark:} Note that there is a misprint in the abstract in  Kr\'ysl
\cite{KryslArchivumRSDR}. Namely, we write there that 
$-\imath l \lambda$ is not a symplectic Killing number instead of 
$\frac{\imath \lambda}{l}$ is not a symplectic Killing spinor number. 
In that paper, $l$ denotes the half of the  dimension of 
the corresponding symplectic manifold.

\bigskip
 
\noindent{\bf Theorem 17:} If $(M,\omega)$ is a symplectic manifold admitting a 
metaplectic structure and  $\nabla$ 
is a Fedosov connection,
then
$$\mbox{kill} = \mbox{Ker} \, T_0 \cap \mbox{Ker}\, \mathfrak{D}.$$ 

\noindent{\it Proof.} See Kr\'ysl \cite{KryslCMUC}.\hfill $\Box$

\bigskip

\noindent{\bf Theorem 18:}
Let $(M^{2n}, \omega)$ be a symplectic manifold admitting a metaplectic 
structure
and $\nabla$ be a Fedosov connection with Ricci tensor $\sigma.$  Let $\phi$ be 
a symplectic Killing spinor field to the symplectic
Killing spinor number $\mu$. Then in a local symplectic frame 
$(U,(e_i)_{i=1}^{2n}),$ we have
$$\Theta^{\sigma} \phi = 2 \mu^2 n \phi.$$

\noindent{\it Proof.} See Kr\'ysl \cite{KryslCMUC}. \hfill $\Box$

\bigskip

As a consequence of this theorem, we have

\bigskip

\noindent{\bf Theorem 19:}
Let $(M,\omega)$ be a symplectic manifold admitting a metaplectic structure and 
$\nabla$ be a Ricci-flat Fedosov connection. Then $kill = \{0\}$ and any symplectic Killing spinor field on $M$ is locally 
covariantly constant.

\noindent{\it Proof.} See Kr\'ysl \cite{KryslCMUC}. 
\hfill $\Box$

\bigskip

\noindent{\bf Remark:} By a locally covariantly constant field $\phi$, we mean $\nabla^S\phi=0$ which implies that
$\phi$ is locally constant if the Kostant's  bundle is trivial. 

\bigskip

\noindent{\bf Theorem 20:} 
Let $n>1,$ $(M^{2n},\omega)$ be a symplectic manifold  admitting a 
metaplectic structure and $\nabla$  be a flat Fedosov
 connection. Then
 \begin{itemize} 
\item[(1)] If $\mu \in  {\it spec}(\mathfrak{D})\setminus 
 \frac{\imath n}{2} \it{kill}$, then 
$\frac{n-1}{n} \mu \in {\it spec}(\mathfrak{R})$.
\item[(2)] If $\phi \in  \mbox{eigen}(\mathfrak{D})\setminus \mbox{kill}
$, then 
$T_0\phi \in \mbox{eigen}(\mathfrak{R}).$
\end{itemize}

\noindent{\it Proof.} See Kr\'ysl \cite{KryslArchivumRSDR}. 
\hfill $\Box$

\bigskip

\noindent{\bf Remark:} For any $\lambda \in \mathbb{C},$  
$\lambda \it{kill}$ denotes the number set $\{\lambda  \alpha, \, \alpha \in \it{kill} 
\}.$


\subsection{First order invariant operators in projective contact geometry} 

Some of the results described above can be modified  to get information for contact projective manifolds which are
more complicated objects to handle than the symplectic ones.
Contact manifolds are models for time-dependent Hamiltonian mechanics. The adjective `projective' is related to 
the fact that we want to deal with unparametrized geodesics rather than with the ones with a fixed parametrization.
Connections that we consider are partial in the sense that they  act on sections of the contact bundle only. 

\bigskip

\noindent{\bf Definition 19:} A {\bf contact manifold} is a manifold $M$ together with a corank one subbundle $HM$ (contact bundle) of the tangent bundle
$TM$ which is not integrable in the Frobenius sense in any point of the manifold, i.e., for each $m\in M,$ there are $\eta_m, \zeta_m \in H_mM$ such that
$[\eta_m, \zeta_m] \notin HM.$ 

\bigskip

Equivalently, $HM$ is a contact bundle if and only if the Levi 
bracket  
$$L(X,Y) =q([X,Y])$$
is non-degenerate. Here $X,Y \in \Gamma(HM)$ and $q: TM \to QM=TM/HM$ denotes the quotient projection onto $QM.$ The Levi bracket induces a tensor field which we denote by the same letter
$L:\bigwedge^2 HM \to QM.$

\bigskip

\noindent{\bf Definition 20:} For a contact manifold $(M,HM),$ 
a partial connection $\nabla: \Gamma(HM) \times \Gamma(HM) \to \Gamma(HM)$ 
is called a {\bf contact connection} if the associated exterior covariant derivative $d^{\nabla}$ on
$\Gamma(\bigwedge^2 HM)$ preserves the kernel of the Levi form, i.e.,
$d^{\nabla}_{\zeta}(\mbox{Ker}\, L) \subseteq \mbox{Ker} \, L$ for any $\zeta \in HM.$ The set of contact connections is denoted by
$\mathcal{C}_M.$
A {\bf contact projective manifold} is a contact manifold $(M, HM)$ together with a set $S_M$ of contact
connections for which the following holds. If $\nabla^1, \nabla^2 \in S_M,$ there exists a differential one-form $\Upsilon \in \Gamma(HM^*)$ such that for any $X,Y \in \Gamma(HM)$
$$\nabla^1_XY - \nabla^2_XY = \Upsilon(X)Y+\Upsilon(Y)X + \Upsilon^{\sharp}(L(X,Y))$$
where $\Upsilon^{\sharp}: QM \to HM$ is a bundle morphism defined by $L(\Upsilon^{\sharp}(\eta),\zeta)=\Upsilon(\zeta)\eta,$ $\zeta \in QM$ and $\eta \in HM.$
Morphisms between  contact projective manifolds $(M,HM, S_M)$ and $(N, HN, S_N)$ are local diffeomorphisms $f: M \to N$ such that $f_*(HM) = HN,$
and for any $\nabla \in S_N,$ the pull-back connection $f^*\nabla \in S_M.$

\bigskip

\noindent{\bf Remark:}   For a contact projective manifold $(M,HM, S_M),$  it is easy to see that the relation $R=S_M \times S_M \subseteq \mathcal{C}_M \times \mathcal{C}_M$ 
on the set of  contact connections $\mathcal{C}_M$ is an equivalence.

\bigskip

Let $(V,\omega)$ be a real symplectic vector space of dimension $2n+2$ and $(e_i)_{i=1}^{2n+2}$ be a symplectic basis.
The action of the symplectic group $G'$ of $(V,\omega)$ on the projectivization
of $V$ is transitive and its stabilizer $P'$ is a parabolic subgroup of $G'$.
We denote the preimages of $G'$ and $P'$ by the covering $\lambda': Mp(2n+2,\mathbb{R}) \to Sp(2n+2,\mathbb{R})$ by
$\widetilde{G}'$ and $\widetilde{P}',$ respectively.

\bigskip

\noindent{\bf Definition 24:}
A {\bf projective contact Cartan geometry} is a Cartan geometry $(\mathcal{G}', \vartheta)$ whose model
is the Klein geometry $G' \to G'/P'$ with $G'$ and $P'$ as introduced above.
We say that a Cartan geometry  is a {\bf metaplectic projective contact Cartan geometry}
if it is modeled on the Klein geometry $\widetilde{G}'/\widetilde{P}'.$ 

\bigskip

\noindent{\bf Remark:} For Cartan geometries, see Sharpe \cite{Sharpe} and \v{C}ap, Slov\'ak \cite{CapSlovak}.
In  \v{C}ap, Slov\'ak \cite{CapSlovak}, a theorem is proved on an
equivalence of the category of the so-called regular normal projective contact 
Cartan geometries and the category
of regular normal projective contact manifolds. See \v{C}ap, Slov\'ak \cite{CapSlovak}, pp. 277 and 410.
See also Fox \cite{Fox}.


\bigskip

The Levi part $\widetilde{G_0}$ of $\widetilde{P}'$ is isomorphic $Mp(2n,\mathbb{R}) \times \mathbb{R}^{\times}$ with the semisimple
part $\widetilde{G_0^{ss}}\simeq \widetilde{G} =  Mp(2n,\mathbb{R})$ and the center isomorphic to the multiplicative group $\mathbb{R}^{\times}.$ 
The Lie algebra $\mathfrak{p}'$ of $\widetilde{P}'$ is graded,
$\mathfrak{p}'  = (\mathfrak{sp}(2n, \mathbb{R}) \oplus \mathbb{R}) \oplus \mathbb{R}^{2n} \oplus  \mathbb{R}$
with  $\mathfrak{g}_0 \simeq \mathfrak{sp}(2n, \mathbb{R}) \oplus \mathbb{R},$ $\mathfrak{g}_1 \simeq \mathbb{R}^{2n} $ and $\mathfrak{g}_2 \simeq \mathbb{R}.$
We denote the Lie algebra of $\widetilde{G}'$ by $\mathfrak{g}'$ and identify it with the Lie algebra $\mathfrak{sp}(2n+2, \mathbb{R}).$
The semi-simple part $\mathfrak{g}_0^{ss}$ of $\mathfrak{g}_0$ is isomorphic $\mathfrak{sp}(2n, \mathbb{R}).$ We denote it by $\mathfrak{g}$ in order 
to be consistent
with the preceding sections.
The grading of $\mathfrak{g}'=\bigoplus_{i=-2}^2 \mathfrak{g}_i,$ $\mathfrak{g}_{-2} \simeq \mathfrak{g}_2$ and $\mathfrak{g}_{-1} \simeq \mathfrak{g}_1,$ can be visualized with respect to the basis $(e_i)_{i=1}^{2n+2}$ by the following block diagonal matrix  of type $(1,n,1) \times (1,n,1)$
$$\mathfrak{g}=\left(
\begin{tabular}{c|ccc|c}
$\mathfrak{g}_{0}$& &$\mathfrak{g}_{1}$&  & $\mathfrak{g}_2$ \\
\hline
&&&&\\
$\mathfrak{g}_{-1}$& &$\mathfrak{g}_0$& &$\mathfrak{g}_1$\\
&&&&\\
\hline $\mathfrak{g}_{-2}$& &$\mathfrak{g}_{-1}$& &$\mathfrak{g}_{0}$\\
\end{tabular}\right).$$
The center of the Lie algebra $\mathfrak{g}_0$ is generated by  
$$Gr=\left(
\begin{tabular}{c|ccc|c}
1& & 0 &  & 0 \\
\hline
&&&&\\
0& &$0$& &0\\
&&&&\\
\hline 0& & 0& &-1\\
\end{tabular}\right)$$
which is usually called the {\bf grading element} because of the 
property $[Gr,X]=jX$ for each $X\in \mathfrak{g}_j$ and $j=-2, \ldots, 2.$



\bigskip

Let $\kappa: (\mathfrak{g}^{\mathbb{C}})^{*} \times (\mathfrak{g}^{\mathbb{C}})^{*} \to \mathbb{C}$ be the dual form to the Killing form of $\mathfrak{g}^{\mathbb{C}} =\mathfrak{sp}(2n,\mathbb{C}).$ 
We choose a Cartan subalgebra $\mathfrak{h}$ of $\mathfrak{g}^{\mathbb{C}}$ and a set of positive roots obtaining the set of fundamental weights
$\{\varpi_i\}_{i=1}^n$ for $\mathfrak{g}^{\mathbb{C}}.$ Further, we set $\langle X,Y \rangle = (4n+4)\kappa(X,Y),$ $X,Y\in (\mathfrak{g}^{\mathbb{C}})^*,$
and define  $$c_{\lambda\nu}^{\mu}=\frac{1}{2}[\langle\lambda,\lambda+2\delta\rangle +\langle\nu,\nu+2\delta\rangle-\langle \mu,\mu+2\delta\rangle]$$ 
for any $\lambda, \mu, \nu \in \mathfrak{h}^*,$ where
$\delta$ is the sum of fundamental weights, or equivalently, the half-sum of positive roots.
For any $\mu \in \mathfrak{h}^*,$ we set 
$$A=\{\sum_{i=1}^{n}\lambda_i\varpi_i|\, \lambda_i \in \mathbb{N}_0, i=1,\ldots, n-1, \lambda_n + 2\lambda_{n-1} +3 >0, \lambda_n \in \mathbb{Z}+\frac{1}{2} \} \subseteq \mathfrak{h}^* \mbox{ and} $$
$$A_{\mu} = A\cap \{\mu + \nu|\, \nu = \pm \epsilon_i, i=1,\ldots, n\}$$
 where $\epsilon_1 =\varpi_1,$ $\epsilon_i = \varpi_i - \varpi_{i-1}, i=2,\ldots, n.$

\bigskip

Considering $\mathbb{C}^{2n}$ with the defining representation of $\mathfrak{g}^{\mathbb{C}} = \mathfrak{sp}(2n,\mathbb{C})$, i.e.,
$\mathbb{C}^{2n} = L(\varpi_1),$ we have the following decomposition.

\bigskip

\noindent{\bf Theorem 21:} For any $\mu \in A,$ the following decomposition into irreducible $\mathfrak{g}^{\mathbb{C}}$-modules 
$$L(\mu) \otimes \mathbb{C}^{2n} = \bigoplus_{\lambda \in A_{\mu}} L(\lambda)$$ holds.

\noindent{\it Proof.} See Kr\'ysl \cite{KryslJOLT1}. 
\hfill $\Box$

\bigskip

\noindent{\bf Remark:} The above decomposition has the same form when we consider the algebra $\mathfrak{sp}(2n,\mathbb{R})$ instead 
of $\mathfrak{sp}(2n,\mathbb{C}).$

\bigskip

The  set $\{L(\lambda)|\, \lambda \in  A\}$ coincides with the set of all infinite dimensional  irreducible
$\mathfrak{g}^{\mathbb{C}}$-modules with  bounded
multiplicities, i.e.,  those irreducible $\mathfrak{sp}(2n,\mathbb{C})$-modules $W$ for which there exists a bound $l\in \mathbb{N}$
 such that for any weight $\nu,$  $\mbox{dim}\, W_{\nu} \leq l.$\footnote{By $W_{\nu}$ we mean the wight space $W_{\nu}=\{w\in W|\, H\cdot w = \nu(H)w \, \mbox{for any } H \in \mathfrak{h}\}.$} See Britten, Hooper, Lemire \cite{BrittenHooperLemire} and Britten, Lemire \cite{BrittenLemire}.

\bigskip

In the next four steps, we define $\widetilde{P}$-modules ${\bf L}(\lambda,c,\gamma)$ for any $\lambda \in A,$ $c\in \mathbb{C}$ and $\gamma \in \mathbb{Z}_2.$

\bigskip

\begin{itemize}
\item[1)]
Let $\mathbb{S}$ and $\mathbb{S}_+$ be the $\mathfrak{g}^{\mathbb{C}}$-modules of 
smooth $\widetilde{K}$-finite vectors of the $Mp(2n, \mathbb{R})$-modules  $L^2(\mathbb{R}^n)$ and $L^2(\mathbb{R}^n)_+$, respectively. Recall that
$L^2(\mathbb{R}^n)$ denotes the Segal--Shale--Weil module and $L^2(\mathbb{R}^n)_+$ is the submodule of even  functions in $L^2(\mathbb{R}^n).$
For any $\lambda \in A,$ there is an irreducible finite dimensional 
$\mathfrak{g}^{\mathbb{C}}$-module $F(\nu)$ with highest weight $\nu \in \mathfrak{h}^*$
such that $L(\lambda)$ is an irreducible summand  in $\mathbb{S}_+ \otimes F(\nu) = \bigoplus_{i=1}^{k} \mathbb{S}_i.$ 
For it, see Britten,  Lemire \cite{BrittenLemire}. Otherwise said, there exists a $j\in \{1,\ldots,k\}$ such that
$L(\lambda) \simeq \mathbb{S}_j.$ The tensor product of the smooth globalization $S=S(\mathbb{R}^n)$ of $\mathbb{S}$ with $F(\nu)$
 decomposes into a finite number of irreducible $\widetilde{G}$-submodules in the corresponding way
$$S_+\otimes F(\nu) = \bigoplus_{i=1}^k S_i$$ i.e., $\mathbb{S}_i$ is the $\mathfrak{g}^{\mathbb{C}}$-module of smooth
$\widetilde{K}$-finite vectors in $S_i.$ We set
${\bf L}(\lambda) = S_j,$ obtaining a $\widetilde{G}$-module.

\item[2)] We let the element $\mbox{exp}(Gr) \in \widetilde{G_0}$ act  by the scalar $\mbox{exp}(c)$ (the conformal weight) on ${\bf L}(\lambda)$ and denote the resulting
structure by ${\bf L}(\lambda, c).$  

\item[3)] Let us consider the element $(1,-1) \in Sp(2n, \mathbb{R}) \times \mathbb{R}^{\times} \subseteq \lambda'(\widetilde{G_0}) \subseteq P$
and  the  preimage  $\Gamma={\lambda'}^{-1}((1,-1)) \subseteq \widetilde{G_0} \simeq Mp(2n,\mathbb{R}) \times \mathbb{R}^{\times}.$
Let us suppose that the element in $\Gamma$ the first component 
of which is the neutral element $e \in Mp(2n,\mathbb{R})$ acts by $\gamma \in \mathbb{Z}_2$ on ${\bf L}(\lambda, c).$ 
 
\item[4)] Finally, the preimage ${\lambda'}^{-1}(G_+) \subseteq \widetilde{P} $ of the unipotent part $G_+$ of $P$ is supposed to act by the identity on 
${\bf L}(\lambda,c).$ We denote the resulting admissible  $\widetilde{P}$-module by 
${\bf L}(\lambda, c, \gamma).$   (See Vogan \cite{Vogan} for the admissibility condition.)
 
\end{itemize}
\bigskip
 

For details on notions in the next definition, see  Slov\'ak, Sou\v{c}ek \cite{SlovakSoucek}.

\bigskip

\noindent{\bf Definition 25:} Let $\mathfrak{G}=(\mathcal{G} \to M,\vartheta)$ be 
a Cartan geometry of type $(G,H)$ and $\mathcal{E}, \mathcal{F} \to M$ be  
 vector bundles associated to the principal $H$-bundle $\mathcal{G} \to M.$
We call a vector space homomorphism $D: \Gamma(\mathcal{E}) \to \Gamma(\mathcal{F})$ a 
{\bf first order invariant differential operator} if
there is a bundle homomorphism $\Phi: J^1 \mathcal{E} \to \mathcal{F}$ such that $Ds = \Phi (s, \nabla^{\vartheta}s)$ for any section $s\in \Gamma(\mathcal{E}),$
 where
$J^1 \mathcal{E}$ denotes the first jet prolongation of $\mathcal{E} \to M$ and $\nabla^{\vartheta}$ is the invariant derivative for
$\mathfrak{G}.$ 

\bigskip

 It is convenient to divide the vector space of first order invariant differential operators by those bundle homomorphisms between $J^1\mathcal{E}$ and $\mathcal{F}$ which act trivially on the tangent space  part of $J^1\mathcal{E}.$ 
We call the resulting vector space the space of {\bf first order invariant operators up to the zeroth order} and denote it 
by $\mbox{Diff}_{\mathfrak{G}}^1(\mathcal{E},\mathcal{F}).$

\bigskip

\noindent{\bf Remark:} Between any bundles induced by  irreducible bounded multiplicities representations introduced above, 
there is at most one such an invariant operator up to a  multiple and up to the operators of  zeroth order.
An equivalent condition for its existence is given in the next theorem. The author obtained it at the infinitesimal level
when  writing his dissertation thesis already. See \cite{KryslDiss}. 

\bigskip

\noindent{\bf Theorem 22:} Let $(\mathcal{G} \to M^{2n+1}, \vartheta)$ be a metaplectic contact projective Cartan geometry,
 $(\lambda, c, \gamma),$ $(\mu, d, \gamma') \in A \times \mathbb{C} \times \mathbb{Z}_2,$ and
   $\mathcal{E} = \mathcal{G} \times_{\widetilde P} {\bf L}(\lambda, c, 
\gamma)$ and 
$\mathcal{F} = \mathcal{G} \times_{\widetilde P} {\bf L}(\mu, d, \gamma')$ be 
the corresponding vector bundles over $M.$
  Then the space 
 $$\mbox{Diff}_{(\mathcal{G} \to M, \vartheta)}^1(\mathcal{E}, \mathcal{F})  
\simeq \left\{\begin{array}{l}
                                                       \mathbb{C}  \quad \hbox{if } \mu \in
                                                       A_{\lambda},  c =  d - 1 
= c_{\lambda\varpi_1}^{\mu} \mbox{ and }
                                                       \gamma = \gamma'
                                                                                                             \\
                                                       0  \quad \hbox{ in other cases.}
                                                       \end{array}
                                                       \right.
                                                       $$

\noindent{\it Proof.} See Kr\'ysl \cite{KryslDGA1}. 
\hfill $\Box$


\subsection{Hodge theory over $C^*$-algebras}

An additive category is called {\bf dagger} if it is equipped with a contravariant 
functor $\ast$ which is the identity on the objects, it is involutive on morphisms, 
$\ast \ast F=F$, and it preserves the identity morphisms, i.e., 
$\ast\mbox{Id}_C = \mbox{Id}_C$ for any object $C$. No compatibility with the additive structure is demanded. See Brinkmann,
Puppe \cite{BrinkmannPuppe}. For a morphism $F,$ we denote $\ast F$ by $F^*.$ For any additive category $\mathcal{C},$ we denote the category of its complexes
by $\mathfrak{K}(\mathcal{C}).$ If $\mathcal{C}$ is an additive and dagger category and $d^{\bullet}=(U^i, d_i)_{i \in \mathbb{Z}} \in \mathfrak{K}(\mathcal{C})$, we set
$\Delta_i = d_i^*d_i+d_{i-1}d_{i-1}^*,$ $i\in \mathbb{Z},$   
and call it the $i$-th {\bf Laplace operator}.

\bigskip

\noindent{\bf Definition 26:} Let $\mathcal{C}$ be an additive and dagger category. We call a complex $d^{\bullet}=(U^i, d_i)_{i\in \mathbb{Z}} \in 
\mathfrak{K}(\mathcal{C})$ of {\bf Hodge-type} if for each $i\in \mathbb{Z}$ 
$$U^i = \mbox{Ker}\,\Delta_i \oplus \mbox{Im} \, d_{i-1} \oplus \mbox{Im}\, d_i^*.$$
We call $d^{\bullet}$ {\bf self-adjoint parametrix possessing}
if  for each $i,$ there exist morphisms $G_i:U^i\to U^i$ and $P_i:U^i\to U^i$ such that
$\mbox{Id}_{U^i} =G_i \Delta_i+ P_i,$ $\mbox{Id}_{U^i} =\Delta_i G_i+ P_i,$ $\Delta_iP_i = 0$ and $P_i = P_i^*.$

\bigskip

\noindent{\bf Remark:}  In the preceding definition, 
we suppose that   the images of the chain maps, the images of their adjoints, and the kernels of the Laplacians
exist as objects in the additive and dagger category $\mathcal{C}.$ The sign $\oplus$ denotes the biproduct in $\mathcal{C}.$
See Weibel \cite{Weibel}, p. 425.
\bigskip

 The first two equations from the definition of a self-adjoint parametrix possessing complex are called the parametrix equations.
Morphisms $P_i$ from the above definition are idempotent as can be seen by composing the first equation with $P_i$ from the right and using the equation
$\Delta_iP_i=0.$ In particular, they are projections. The operators $G_i$ are called the Green operators. 

\bigskip

\noindent{\bf Definition 27:}
Let $(A,*_A,|\,|_A) $ be a $C^*$-algebra and $A^+$ be the positive cone of $A,$ i.e., the set of all hermitian elements ($*_A a=a$) in $A$ whose spectrum  is contained in  the non-negative real numbers.  
 A tuple $(U,(,))$ is called a {\bf pre-Hilbert $A$-module} if $U$ is a right module over the complex associative algebra $A,$ and
$(,):U \times U \to A$ is an $A$-sesquilinear map such that for all $u,v\in U,$ $(u,v)=*_A(v,u),$ $(u,u)\in A^+,$ and $(u,u)=0$ implies $u=0.$
A pre-Hilbert module is called a {\bf Hilbert $A$-module} if it is complete with respect to the norm $|u|=\sqrt{|(u,u)|_A},$ $u \in U.$
A pre-Hilbert $A$-module morphism between $(U,(,)_U)$ and $(V,(,)_V)$ is any continuous $A$-linear map $F: U \to V.$

\bigskip

\noindent{\bf Remark:} We consider that $(,)$ is antilinear in the left variable and linear in the right one as it is usual in physics.

\bigskip

An adjoint of a morphism $F:U\to V$ acting between pre-Hilbert modules $(U,(,)_U)$ and $(V,(,)_V)$ is a morphism $F^*:V \to U$ that satisfies the condition 
$(Fu,v)_V= (u, F^*v)_U$ for any $u \in U$ and $v\in V.$
The category of pre-Hilbert and Hilbert $C^*$-modules and adjointable morphisms is an additive and dagger category. The dagger functor is the adjoint 
on morphisms. For any $C^*$-algebra $A,$ we denote the categories of pre-Hilbert $A$-modules and Hilbert $A$-modules and adjointable morphisms by $PH_A^*$ and
$H_A^*,$ respectively. In both of these cases, the dagger structure is compatible with the additive structure.

\bigskip

To any complex $d^{\bullet}=(U^i, d_i)_{i\in \mathbb{Z}} \in \mathfrak{K}(PH_A^*),$ the cohomology groups $H^i(d^{\bullet})= \mbox{Ker}\, d_i / \mbox{Im}\, d_{i-1}$ are assigned which
are $A$-modules and which we consider to be equipped with the canonical quotient topology. They are pre-Hilbert $A$-modules with respect to the restriction 
of $(,)_{U_i}$ to $\mbox{Ker}\, d_i $ if and only if $\mbox{Im}\, d_{i-1}$ has 
 an $A$-orthogonal complement in $\mbox{Ker}\, d_i.$ 

\bigskip

We have the following

\bigskip

\noindent{\bf Theorem 23:}
Let $d^{\bullet} = (U^i, d_i)_{i\in \mathbb{Z}}$ be a self-adjoint 
parametrix possessing complex in $PH_A^*.$ Then for any $i \in \mathbb{Z}$
\begin{itemize}
\item[1)]$d^{\bullet}$ is of Hodge-type 
\item[2)]$H^i(d^{\bullet})$ is isomorphic to $\mbox{Ker}\, \Delta_i \mbox{ as a pre-Hilbert $A$-module}$
\item[3)]$\mbox{Ker} \, d_i =\mbox{Ker} \, \Delta_i \oplus \mbox{Im} \, d_{i-1}$
\item[4)]$\mbox{Ker} \, d_i^* = \mbox{Ker} \, \Delta_{i+1} \oplus \mbox{Im} \, d_{i+1}^*$
\item[5)]$\textrm{Im}\, \Delta_i  = \textrm{Im} \, d_{i-1} \oplus \textrm{Im} \,  d_i^*.$
\end{itemize}

\noindent{\it Proof.} See Kr\'ysl \cite{KryslAGAG2}. 
\hfill $\Box$

\bigskip

\noindent{\bf Remark:}
If the image of $d_{i-1}$ is not closed, the quotient topology on the cohomology group $H^i(d^{\bullet})$ is non-Hausdorff and in particular,
it is not in $PH_A^*.$ See, e.g., von Neumann \cite{vonNeumann} on the relevance of topology for state spaces. See also Kr\'ysl \cite{KryslJGP2} for  
further references and for a relevance of our topological observation (Theorem 23 item 2) to the  basic principles of the  so-called Becchi--Rouet--Stora--Tyutin (BRST) quantization.

\bigskip

\noindent{\bf Theorem 24:} Let $d^{\bullet} = (U^{i}, d_{i})_{i\in \mathbb{Z}}$ be a complex of Hodge-type in 
$H_A^*$, then $d^{\bullet}$ is self-adjoint parametrix possessing.

\noindent{\it Proof.} See Kr\'ysl \cite{KryslJGP2}. \hfill $\Box$

\bigskip

\noindent{\bf Definition 28:}
Let $M$ be a smooth manifold, $A$ be a $C^*$-algebra and $\mathcal{F} \to M$ be a Banach bundle with a 
smooth atlas such that each of its maps targets onto a fixed Hilbert $A$-module (the typical fiber).
If  the transition 
functions of the atlas are Hilbert $A$-module automorphisms, we call $\mathcal{F} 
\to M$ an $A$-Hilbert bundle.  
We call an $A$-Hilbert bundle $\mathcal{F} \to M$  {\bf finitely generated projective}
if the typical fiber is a finitely generated projective Hilbert $A$-module.

\bigskip

For further information on analysis on $C^*$-Hilbert bundles, we refer to Solovyov, Troitsky 
\cite{SolovyovTroitsky}, Troitsky \cite{Troitsky} and Schick \cite{Schick}. In the paper of Troitsky, complexes  are
 treated with an allowance of the so-called `compact' perturbations.


 

\bigskip

\noindent{\bf Theorem 25:}
Let $M$ be a compact manifold, $A$ be a  $C^*$-algebra and $D^{\bullet} = 
(\Gamma(\mathcal{F}^i), 
D_i)_{i \in \mathbb{Z}}$ be 
an elliptic complex on finitely generated projective $A$-Hilbert bundles 
over $M$. Let for each $i\in \mathbb{Z},$
 the image of  $\Delta_i$ be closed in 
$\Gamma(\mathcal{F}^i).$ 
Then for any $i\in \mathbb{Z}$
\begin{itemize}
\item[1)] $D^{\bullet}$  is of Hodge-type 
\item[2)] $H^i(D^{\bullet})$ is a finitely generated projective Hilbert 
$A$-module isomorphic to $\mbox{Ker}\,\Delta_i$ as a Hilbert $A$-module
\item[3)] $\mbox{Ker} \, D_i = \mbox{Ker} \, \triangle_{i} \oplus \mbox{Im} \,
D_{i-1}$
\item[4)] $\mbox{Ker} \, D_i^* = \mbox{Ker} \, \triangle_{i+1} \oplus \mbox{Im} 
\, D_{i+1}^*$
\item[5)] $\mbox{Im}\, \Delta_i = \mbox{Im}\, D_{i-1} \oplus \mbox{Im} \, D_i^*.$
\end{itemize}

\noindent{\it Proof.} See Kr\'ysl \cite{KryslAGAG2}. \hfill $\Box$

\bigskip

Let $H$ be a Hilbert space. Any $C^*$-subalgebra of the 
$C^*$-algebra of  compact operators on $H$ is called a 
{\bf $C^*$-algebra of compact operators}.

\bigskip

For $C^*$-algebras of compact operators, we have the following analogue of the Hodge theory for elliptic complexes of 
operators on finite rank vector bundles over compact manifolds.

\bigskip

\noindent{\bf Theorem 26:} Let $M$ be a compact manifold, $K$ be a  $C^*$-algebra of compact operators and $D^{\bullet} = 
(\Gamma(\mathcal{F}^i), D_i)_{i \in \mathbb{Z}}$ be 
an elliptic complex on finitely generated projective $K$-Hilbert bundles 
over $M$.  
If  $D^{\bullet}$ is elliptic,
then  for each $i\in \mathbb{Z}$
\begin{itemize}
\item[1)]  $D^{\bullet}$ is of Hodge-type
\item[2)] The cohomology group $H^i(D^{\bullet})$ is a finitely generated 
projective Hilbert $K$-module isomorphic to the Hilbert $K$-module $\textrm{Ker} \, \Delta_i.$ 
\item[3)] $\textrm{Ker}\, D_{i} = \textrm{Ker} \, \Delta_i \oplus 
\textrm{Im} \, D_{i-1}$
\item[4)] $\textrm{Ker}\, D_{i}^* = \textrm{Ker} \, \Delta_{i+1} \oplus 
\textrm{Im} \, D_{i+1}^*$
\item[5)] $\textrm{Im}\, \Delta_i  = \textrm{Im} \, D_{i-1} \oplus 
\textrm{Im} \,  D_i^*$
\end{itemize}
 
\noindent{\it Proof.} See \cite{KryslJGP2}. \hfill $\Box$

\bigskip

\noindent{\bf Remark:} In particular, we see that the cohomology groups share  properties of
the fibers.

\newpage



`

\bibliographystyle{amsplain}

\newpage

\section{Selected author's articles}

\bigskip

\begin{itemize}

\item[]  Decomposition of a tensor product of a 
higher symplectic spinor module and the defining representation of 
$\mathfrak{sp}(2n,\mathbb{C})$, Journal of Lie Theory, Vol. 17 (1) (2007), 
63--72. 

\bigskip

\item[] Classification of 1st order symplectic spinor operators in contact projective geometries, Differential Geometry and Application, Vol. 26 (3)
(2008), 553--565.
 
\bigskip

\item[] Structure of the curvature tensor on symplectic spinors, Journal Geometry and Physics, Vol. 60(9) (2010), 1251--1261. 

\bigskip

 \item[]  Complex of twistor operators in spin symplectic geometry, Monatshefte f\"ur Mathematik, Vol. 161 (4) (2010), 381--398.  

\bigskip

\item[] Ellipticity of symplectic twistor 
complexes, Archivum Math., Vol. 44 (4) (2011), 309--327.
  
\bigskip

 \item[]   Symplectic Killing spinors, Comment. Math. 
Univ. Carolin., Vol. 53 (1) (2012) 19--35. 
\bigskip

 \item[] Howe duality for the metaplectic group acting on symplectic spinor valued forms, Journal of Lie theory, Vol. 22 (4) (2012), 1049--1063.
\bigskip


 \item[] Cohomology of the de Rham complex twisted by the oscillatory representation, Diff. Geom. Appl., Vol. 33 (Supplement) (2014), 290--297. 
\bigskip


 \item[] Hodge theory for complexes over 
$C^*$-algebras with an application to $A$-ellipticity, Annals Glob. Anal. Geom., 
Vol. 47 (4) (2015), 359--372. 
\bigskip

\item[] Elliptic complexes over $C^*$-algebras of compact operators, Journal of Geometry and Physics, Vol. 101 (2016), 27--37.

\end{itemize}

\end{document}